\documentclass[a4paper,leqno]{amsart}
\renewcommand{\today}{2012/12/20}
\usepackage[dvips]{graphicx}
\usepackage[usenames]{color}
\newcommand{\vect}[1]{\boldsymbol{#1}}

\newcommand{\imag}{\mathrm{i}}
\newcommand{\D}{\mathbb{D}}
\newcommand{\DD}{\mathcal{D}}
\newcommand{\C}{\mathbb{C}}
\newcommand{\R}{\mathbb{R}}

\newcommand{\M}{\operatorname{M}}
\newcommand{\A}{\mathcal{A}}
\newcommand{\SO}{\operatorname{SO}}
\newcommand{\SL}{\operatorname{SL}}

\newcommand{\U}{\operatorname{U}}
\renewcommand{\O}{\operatorname{O}}
\newcommand{\id}{\operatorname{id}}

\newcommand{\Length}{\operatorname{Length}}
\newcommand{\trans}[1]{{#1}^t}
\renewcommand{\Im}{\operatorname{Im}}
\renewcommand{\Re}{\operatorname{Re}}

\newcommand{\hinner}[2]{\left\langle #1,#2\right\rangle}
\newcommand{\inner}[2]{\left(#1,#2\right)}
\newcommand{\dist}{\operatorname{dist}}
\newcommand{\pmt}[1]{{\begin{pmatrix} #1  \end{pmatrix}}}
\renewcommand{\theenumi}{(\arabic{enumi})}
\renewcommand{\labelenumi}{(\arabic{enumi})}
\numberwithin{equation}{section}
\newtheorem{theorem}{Theorem}[section]
\newtheorem{corollary}[theorem]{Corollary}
\newtheorem{lemma}[theorem]{Lemma}

\newtheorem{proposition}[theorem]{Proposition}

\theoremstyle{definition}

\makeatletter
\def\paragraph{\@startsection{paragraph}{4}%
  \z@\z@{-\fontdimen2\font}%
  {\bfseries}}
\def\subparagraph{\@startsection{paragraph}{4}%
  \z@\z@{-\fontdimen2\font}%
  {\itshape}}
\makeatother
\title[A bounded null curve]{%
   A construction of\\
   a complete bounded null curve in $\C^3$
}
\author{L. Ferrer}
\address[Leonor Ferrer]{
  Departamento de Geometr\'\i{}a y Topolog\'\i{}a,
  Universidad de Granada,
  18071 Granada, Spain.
}
\email{lferrer@ugr.es}
\author{F. Mart\'\i{}n}
\address[Francisco Mart\'\i{}n]{
  Departamento de Geometr\'\i{}a y Topolog\'\i{}a,
  Universidad de Granada,
  18071 Granada, Spain.
}
\email{fmartin@ugr.es}
\author{M. Umehara}
\address[Masaaki Umehara]{%
   Department of Mathematical and Computing Sciences,
   Tokyo Institute of Technology
   2-12-1-W8-34, O-okayama, Meguro-ku,
   Tokyo 152-8552, Japan.
}
\email{umehara@is.titech.ac.jp}
\author{K. Yamada}
\address[Kotaro Yamada]{%
   Department of Mathematics,
   Tokyo Institute of Technology
   2-12-1-H-7, O-okayama, Meguro-ku,
   Tokyo 152-8551, Japan.
}
\email{kotaro@math.titech.ac.jp}
\date{\today}
\begin{document}
\begin{abstract}
We construct a complete bounded immersed null holomorphic
curve in $\C^3$, which is a recovery of the previous
paper of the last three authors on this subject.
\end{abstract}
\maketitle
\section*{Introduction}
The study of global properties of complete complex null-curves is
interesting from different points of view.
Firstly, the real and imaginary part of such a curve are complete
minimal surfaces in $\R^3$. 
Secondly, there exists a close relationship between null curves in
$\C^3$ and constant mean curvature $H=1$ surfaces in hyperbolic 
$3$-space.

An important problem in the global theory of complete null curves is the
so called Calabi-Yau problem,  
which deals with the existence of complete null-curves inside a ball of
$\C^3$.  
This problem was approached firstly in \cite[Theorem A]{MUY}
using similar ideas to those used by Nadirashvili in \cite{Nadi} to
solve the Calabi-Yau conjecture in $\R^3$.
Unfortunately, the paper \cite{MUY} has a mistake, and the first
examples of complete bounded null curves in $\C^3$ where provided using
other approximation, by Alarc\'on and L\'opez  \cite{AL}. 
Very recently, Alarc\'on and Forstneri\v{c} have got the most general
results in this line (see \cite{AF1,AF2}).

The purpose of this paper is to show that similar ideas to those given
in  \cite{MUY} can be used  to produce examples of complete bounded
null holomorphic disks in a ball of $\C^3$:
In \cite{MUY}, Mart\'i{}n, Umehara and Yamada tried to construct a
bounded holomorphic curve in $\SL(2,\C)$ and used this example 
to get the desired bounded disk in $\C^3$. 
However, in this paper, we construct the bounded null curves directly in
$\C^3$. 
In this aspect, our strategy is similar to that used by Alarc\'on and
L\'opez in \cite{AL}.
Although, as we mentioned before, these examples have been generalized
in Alarc\'on and Forstneri\v{c} \cite{AF2} by using different (and
powerful) methods, we think that the arguments and techniques exhibited
in this paper is different from \cite{AL, AF1, AF2}, and
might be of use in the solution of other questions related
to the Calabi-Yau problem in different settings.

As applications of Theorem A in \cite{MUY},
the following objects were constructed;
\begin{enumerate}
\item complete bounded minimal surfaces in  the Euclidean $3$-space $\R^3$
      (\cite[Theorem A]{MUY}),
\item complete bounded holomorphic curves in $\C^2$ 
      (\cite[Corollary~B]{MUY}),
\item weakly complete bounded maximal surfaces in the Lorentz-Minkowski
      3-space $\R_1^3$
      (\cite[Corollary D]{MUY}),
\item complete bounded null curves in $\SL(2,\C)$
      (\cite[Theorem C]{MUY}),
\item complete bounded constant mean curvature
      one surfaces in the hyperbolic 3-space $H^3$
      (\cite[Theorem C]{MUY}).
\end{enumerate}
We also constructed 
higher genus examples of the first three objects
in \cite{MUY2}. 
All of these applications in \cite{MUY} and \cite{MUY2} are correct
as a consequence.

\section{The Main Theorem and the Key Lemma}
\label{sec:main}
We denote by $\inner{~}{~}$ (resp.\ $\hinner{~}{~}$)
the $\C$-bilinear inner product 
(resp.\ the Hermitian inner product) of $\C^3$:
\begin{equation}
\label{eq:inner}
     \inner{\vect{x}}{\vect{y}}:=
     x_1y_1+x_2y_2+x_3y_3,\qquad
     \hinner{\vect{x}}{\vect{y}}:=\inner{\vect{x}}{\overline{\vect{y}}},
\end{equation}
where $\vect{x}=(x_1,x_2,x_3)$, 
$\vect{y}=(y_1,y_2,y_3)\in\C^3,$ and $\overline{\vect{y}}$
denotes the complex conjugate of $\vect{y}$.
We identify an element of $\C^3$ with a column vector when the 
matrix product is used.
The Hermitian norm of $\C^3$ is denoted by
$|\vect{x}|:=\sqrt{\hinner{\vect{x}}{\vect{x}}}$ for $\vect{x}\in\C^3$.
In particular, it holds that 
\begin{equation}\label{eq:schwarz}
 \bigl|\inner{\vect{x}}{\vect{y}}\bigr|=
 \bigl|\hinner{\vect{x}}{\overline{\vect{y}}}\bigr|
 \leq |\vect{x}|\,|\overline{\vect{y}}|=|\vect{x}|\,|\vect{y}|.
\end{equation}
Let 
$\M(3,\C)$ (resp.\ $\M(3,\R)$) be the set of complex (resp.\ real)
$(3\times 3)$-matrices. Moreover, we will use the following notation for the
set of complex (resp.\  special) orthogonal matrices
\begin{align*}
    &\O(3,\C):=\{A\in\M(3,\C)\,;\,\trans{A}A=\id\},\\
    &\biggl(\text{resp.\ }\SO(3)
      :=\{A\in\M(3,\R)\,;\,\trans{A}A=\id,  \,\, 
                           \det A=1\}\biggr),
\end{align*} where $\trans{A}$ means the transposed matrix of $A$.
As usual, we denote $\U(3):=\{A\in\M(3,\C)\,;\,A^*A=\id\}$,
where $A^*$ is the conjugate transposed matrix of $A$.
For each $A\in\M(3,\C)$, we define the {\em matrix norm\/} as
\begin{equation}\label{eq:matrix-norm}
    \Vert{}A\Vert{}:=\sup_{\vect{x}\in\C^3\setminus\{\vect{0}\}}
          \frac{|A\vect{x}|}{|\vect{x}|}.
\end{equation}
 If  $A\in \M(3,\C)$ is a non-singular matrix, 
\begin{equation}
\label{eq:norm-ineq}
  \dfrac{1}{\Vert{}A^{-1}\Vert{}}|\vect{x}|\leq|A\vect{x}|\leq 
  \Vert{}A\Vert{}\,|\vect{x}|
\end{equation}
 holds.
It is well-known that
\begin{equation}\label{eq:matrix-norm-1}
     \Vert{}A\Vert{} = \sqrt{\max\{\mu_1,\mu_2,\mu_3\}}
\qquad (A\in\M(3,\C))
\end{equation}
holds, where $\mu_j \in \R$ $(j=1,2,3)$ are 
 the eigenvalues of positive semi-definite Hermitian 
 matrix $A^*A$.

\bigskip
A holomorphic map $F\colon{}D\to\C^3$, defined on a domain $D\subset\C$, is 
a null immersion if and only if
\begin{equation}\label{eq:phi-f}
   \inner{\varphi_F}{\varphi_F}=0 \quad\text{and}\quad
   |\varphi_F|^2=\hinner{\varphi_F}{\varphi_F}>0,
   \qquad\text{where}\quad
      \varphi=\varphi_F:=\frac{dF}{dz},
\end{equation} where $z$ is the 
canonical complex coordinate of $\C$.
In this case the pull-back of the Hermitian metric of $\C^3$ by $F$ is 
expressed as 
\begin{equation}\label{eq:induced}
   ds^2_F:=\hinner{dF}{dF}=|\varphi_F|^2\,|dz|^2,
\end{equation}
which is called the {\em induced metric\/} of $F$.
For a holomorphic null immersion $F\colon{}D\to\C^3$, the first
equality of \eqref{eq:phi-f} implies that
there exist a meromorphic function $g$ and a holomorphic function $\eta$
such that
\begin{equation}\label{eq:weierstrass}
   \varphi_F 
       = \frac{1}{2}\,\bigl(1-g^2,\imag(1+g^2),2g\bigr)\eta
       \qquad \left(\imag = \sqrt{-1}\right).
\end{equation}
We call $(g,\eta)$ the {\em Weierstrass data\/} of $F$.
Using these data, the induced metric \eqref{eq:induced} is 
expressed as 
\begin{equation}\label{eq:induced-weier}
   ds^2_F = 
          \frac{1}{2}(1+|g|^2)^2\,|\eta|^2\,|dz|^2.
\end{equation}

Throughout this paper, we denote the open (resp.\ closed) disc on $\C$
centered at $0$ with radius $r$ by 
\begin{equation}\label{eq:disc}
  \D_r := \{z\in\C\,;\,|z|<r\},\qquad
  (\text{resp.\ }
  \overline{\D}_r := \{z\in\C\,;\,|z|\leq r\})\qquad (r>0).
\end{equation}
The goal of this paper is to prove the following
\begin{theorem}[The Main Theorem]\label{thm:main}
 There exists a holomorphic null immersion 
 $X\colon{}\D_1\to\C^3$ such that the induced metric
 $ds^2_X$ is complete, and the image $X(\D_1)$
 is bounded in $\C^3$.
\end{theorem}
The main theorem can be proved by the following 
assertion in the same way as \cite{Nadi,CR}:
\begin{proposition}\label{prop:main}
 Let $X\colon{}\overline{\D}_1\to\C^3$ be a holomorphic
 null immersion of the closed disc $\overline{\D}_1\subset\C$
 into $\C^3$.
 Suppose that there exists positive numbers $\rho$ and $r$
 such that
\begingroup
 \renewcommand{\theenumi}{(X-\arabic{enumi})}
 \renewcommand{\labelenumi}{(X-\arabic{enumi})}
 \begin{enumerate}
 \item\label{key:ass:1} 
      $X(0)=0$,
 \item\label{key:ass:2} 
      $(\D_1,ds^2_X)$ contains the geodesic disc
      centered at $0$ and of radius $\rho$,
 \item\label{key:ass:3}
      and $|X|\leq r$ holds on 
      $\overline{\D}_1$.
 \end{enumerate}
 \endgroup
 Then, for an arbitrary given positive numbers $\varepsilon$ and $s$,
 there exists 
 a holomorphic null immersion $Y\colon{}\overline{\D}_1\to\C^3$ satisfying
 \begingroup
 \renewcommand{\theenumi}{(Y-\arabic{enumi})}
 \renewcommand{\labelenumi}{(Y-\arabic{enumi})}
 \begin{enumerate}
  \item\label{key:concl:d1}
       $|\varphi_Y-\varphi_X|<\varepsilon$ and 
       $|Y-X|<\varepsilon$ hold on
       $\D_{1-\varepsilon}$, where 
       $\varphi_X=dX/dz$ and 
       $\varphi_Y=dY/dz$,
  \item\label{key:concl:d2}
       $(\D_1,ds^2_{Y})$ contains the geodesic disc $\DD$
       centered at $0$ with radius  $\rho+s$,
  \item\label{key:concl:d3}
       and on the boundary $\partial\DD$ of the geodesic disc
       $\DD$ in \rm{\ref{key:concl:d2}},
       it holds that 
       $|Y|\leq  \sqrt{r^2+s^2}+\varepsilon$.
 \end{enumerate}
 \endgroup
\end{proposition}
This proposition is a consequence of the following Key Lemma.
(The proof of Proposition~\ref{prop:main} is given in
 Section~\ref{sec:proof-main}.)
To explain it, we define three constants
\[
  N=N(\rho,r,\mu,\nu,s,\varepsilon),\qquad
  C_1=C_1(\rho,r,\mu,\nu,s,\varepsilon)\quad\text{and}\quad
  C_2=C_2(\rho,r,\mu,\nu,s,\varepsilon)
\]
depending on six positive constants
$\rho$, $r$, $\mu$, $\nu$, $s$, $\epsilon$.
Here $\rho$ and $r$  have been already given in
\ref{key:ass:2} and \ref{key:ass:3},
and we will fix $\mu,\nu$ in the statement of Lemma~\ref{lem:key}. 
The remaining two constants $s$, $\epsilon$ are arbitrary in
the statement of Lemma \ref{lem:key}, 
but will coincide with the corresponding constants as in 
Proposition~\ref{prop:main}.

The constants $C_1$ and $C_2$ are set as
\begin{equation}\label{eq:C1C2}
  C_1:= \frac{\nu}{5},
 \qquad
  C_2:= 6(\mu^2 + 2\mu+2).
\end{equation}
Next, we set
\begin{equation}\label{eq:c123}
\begin{aligned}
  &c_1:=6\mu^2+12\mu+8, \quad c_2:= 3\mu+\frac{2\varepsilon(\rho+s)}{C_1},
\\
  &c_3:=
     \frac{%
        s\alpha+\frac{\alpha^2}{2}+2 r\varepsilon+2\varepsilon^2}{%
        \sqrt{r^2+s^2}}
     \quad
     \biggl(\alpha:=
             c_2+5\varepsilon+(r+2\varepsilon)\sqrt{2C_2}\biggr).
\end{aligned}
\end{equation}
We then choose an integer $N$ so that
it satisfies the following four inequalities;
\begin{align}
\label{eq:N1}
 N&\geq \max\left\{
              36,~\frac{2\varepsilon}{\nu},~\varepsilon,~
              (12\mu)^2,~
              \left[
              \frac{%
                2^5(3\mu+\varepsilon)}{\nu}\left(2+\frac{6\mu+\varepsilon}{3\nu}
                 \right)
              \right]^4
            \right\}, \\
\label{eq:N3}
 N&\geq \max\left\{
               \frac{3}{\varepsilon},~
               \left(\frac{2\varepsilon}{C_1}\right)^4\!\!\!,
               ~~
               \left(
                \frac{1}{\nu}
                \left(\varepsilon+\frac{C_1}{2}\right)
               \right)^{4/3}\!\!\!,~~
              \left(\frac{2(\rho+s)}{C_1}\right)^4\right\},\\
\label{eq:N4}
 N&\geq \max\left\{
           \left(\frac{c_3+2\varepsilon}{\varepsilon}\right)^4
                  \!\!\!,~~
           \left(\frac{1+c_2+6\mu+3\varepsilon}{\varepsilon}\right)^4
        \right\}.
\end{align}

\begin{lemma}[The Key Lemma]\label{lem:key}
 Assume a holomorphic null immersion 
 $X\colon{}\overline{\D}_1\to\C^3$ and positive real numbers
 $\rho$ and $r$ satisfy {\rm\ref{key:ass:1}--\ref{key:ass:3}}.
 We set 
 \begin{equation}\label{eq:constants}
   \nu:=\min_{\overline{\D}_1}|\varphi_X|>0,\qquad
   \mu:=\max\left\{1+\max_{\overline{\D}_1}|\varphi_X|,~
	   \max_{\overline{\D}_1}|\varphi_X'|\right\},
 \end{equation} 
 where
 $\varphi_X:=X'=dX/dz$ and 
 $\varphi_X':=d\varphi_X/dz$.
 For an arbitrary positive number $\varepsilon$ and $s$,
 we take positive constant $C_1$, $C_2$ and positive integer $N$
 as in \eqref{eq:C1C2}, \eqref{eq:N1}--\eqref{eq:N4}.
 Then there exist a sequence $\{F_j\}_{j=0,\dots,2N}$ of 
 holomorphic null immersions 
 $F_j\colon{}\overline{\D}_1\rightarrow\C^3$
 and a sequence $\{\vect{v}_j\}_{j=1,\dots,2N}$
 of unit vectors in $\C^3$  which satisfy the following assertions,
 where 
 the compact set $\omega_j\subset\C$, 
 an open neighborhood $\varpi_j$ of $\omega_j$
 and the ``base point'' $\zeta_j$ of $\varpi_j$
 are as in \eqref{eq:omega} and \eqref{eq:zeta}
 in Appendix~\ref{app:labyrinth}, and
 \begin{equation}\label{eq:phi}
  \varphi_l = \frac{dF_l}{dz}
   \qquad (l=0,\dots,2N)
 \end{equation}
\begingroup
 \renewcommand{\theenumi}{(K-\arabic{enumi})}
 \renewcommand{\labelenumi}{(K-\arabic{enumi})}
 \begin{enumerate}
 \setcounter{enumi}{-1}
  \item\label{key:concl:0}
       $F_0=X$.
  \item\label{key:concl:1}
       $F_l(0)=0$ $(l=0,\dots,2N)$.
  \item\label{key:concl:2}
       $|\varphi_l-\varphi_{l-1}|\leq \frac{\varepsilon}{2N^{2}}$
holds on $\overline{\D}_1\setminus\varpi_l$ for each
$l=1,\dots,2N$. 
  \item\label{key:concl:3}
       The inequality 
       \[
	    |\varphi_l|\geq \begin{cases}
			  C_1 N^{9/4}\qquad &\text{on $\omega_l$},\\
			  C_1 N^{-3/4}\qquad &\text{on $\overline{\varpi}_l$}
			 \end{cases}
	\]
       holds for each $l=1,\dots,2N$.
  \item\label{key:concl:4}
       $|\inner{\vect{v}_l}{\vect{v}_l}|\geq {1}/{{N}^{1/4}}$
       for each $l=1,\dots,2N$.
  \item\label{key:concl:5}
       $|F_{l-1}(\zeta_l)|<1/\sqrt{N}$, or
       \[
          \left|\hinner{
                  \frac{F_{l-1}(p)}{|F_{l-1}(p)|}
                }{\vect{v}_l}\right|
                \geq 1-\frac{C_2}{\sqrt{N}}
           \qquad (\text{on $\overline{\varpi}_l$})
       \]
       holds for each $l=1,\dots,2N$.
  \item\label{key:concl:6}
       $\hinner{F_{l-1}}{\vect{v}_l} =  \hinner{F_l}{\vect{v}_l}$
       holds on $\overline{\D}_1$ for each $l=1,\dots,2N$.
 \end{enumerate}
 \endgroup%
\end{lemma}
In the proof of the Key Lemma~\ref{lem:key}, we use the notion of
Gauss maps of holomorphic null immersions:
Let $F\colon{}D\to\C^3$  be a holomorphic null immersion.
Then both the real part $\Re F$ and the imaginary part $\Im F$
give conformal minimal immersions into $\R^3$ with the same Gauss 
map. So we call the Gauss map  $G\colon{}D\to S^2$ 
of both $\Re F$ and $\Im F$ the {\em Gauss map\/} of $F$,
where $S^2\subset \R^3$ is the unit sphere.
Then $G$ is expressed as
\begin{equation}\label{eq:gauss}
     G = \frac{-\imag\, (\varphi\times\bar\varphi)}{|\varphi|^2}\colon{}
      D\longrightarrow S^2\subset \R^3
    \qquad\left(\varphi=\frac{dF}{dz}\right),
 \end{equation}
because \eqref{eq:phi-f} implies that
$|\varphi\times\bar\varphi|=|\varphi|^2$,
where  ``$\times$'' denotes the complexification of the vector
product of  $\R^3$.
Using the Weierstrass data \eqref{eq:weierstrass},
$G$ is expressed as 
\begin{equation}\label{eq:gauss-weier}
   G = \left(
               \frac{2\Re g}{1+|g|^2},
               \frac{2\Im g}{1+|g|^2},
               \frac{|g|^2-1}{1+|g|^2}
              \right).
\end{equation}
That is, $g=\pi_S\circ G$,
where $\pi_S\colon{}S^2\to\C\cup\{\infty\}$ is the stereographic
projection from the north pole.
\section{Preliminary estimates}
\label{sec:prelim-2}
Let $F_0=X\colon{}\overline{\D}_1\to \C^3$ be a holomorphic null immersion
as in the assumption of the Key Lemma~\ref{lem:key}.
Here, we prepare some basic properties of $\{F_j\}_{j=0,\dots,2N}$ 
in the conclusion of the Key Lemma~\ref{lem:key}.
\begin{lemma}\label{lem:fj-est}
 If   {\rm\ref{key:concl:1}} and {\rm \ref{key:concl:2}}  
in the Key Lemma~\ref{lem:key} are satisfied for $l\in \{1,\ldots,2N\}$ then 
 \[
     |F_l-F_{l-1}|\leq \frac{\varepsilon}{N^2}
     \qquad\text{on $\overline{\D}_1\setminus\varpi_l$}.
 \]
\end{lemma}
\begin{proof}
 Let $p\in \overline{\D}_1\setminus\varpi_l$.
 Then there exists a path $\gamma$ in $\overline{\D}_1\setminus\varpi_l$
 joining $0$ and $p$ whose Euclidean length is not greater 
 than $1+\frac{\pi}{N}$ (see Lemma~\ref{lem:path} in 
 Appendix~\ref{app:labyrinth}).
 Thus, we have
 \begin{alignat*}{2}
    |F_{l}(p)-F_{l-1}(p)|&=
    \left|\int_{\gamma} \bigl(\varphi_l(z)-\varphi_{l-1}(z)\bigr)
            \,dz\right|\qquad &&\text{(by \ref{key:concl:1})}\\
     &\leq 
             \int_{\gamma}
                \bigl|\varphi_l(z)-\varphi_{l-1}(z)\bigr|\,|dz|
     \leq \Length_{\C}(\gamma) \frac{\varepsilon}{2N^2}\quad
             &&\text{(by \ref{key:concl:2})}\\
     &\leq \left(1+\frac{\pi}{N}\right)\frac{\varepsilon}{2N^2}
      \leq 2\cdot\frac{\varepsilon}{2N^2}=\frac{\varepsilon}{N^2}
           \qquad&&\text{(by \eqref{eq:N1})},
 \end{alignat*}
 where $\Length_{\C}(\gamma)$ is the length of
 $\gamma$ with respect to the metric $|dz|^2$ on $\C$.
\end{proof}
\begin{lemma}\label{lem:bound-j}
 Fix an integer $j$ {\rm(}$1\le j\le 2N${\rm)}.
 If $F_0$, $F_1$, \dots $F_{j-1}$ satisfy
 {\rm \ref{key:concl:0}} and {\rm \ref{key:concl:2}} of the Key
 Lemma~\ref{lem:key}.
 Then 
 \[
    |\varphi_{j-1}|\geq \frac{\nu}{2}\qquad\text{and}\qquad
    |\varphi_{j-1}|\leq \mu
 \qquad
  \left(
  \text{on 
    $\overline{\D}_1\setminus(\varpi_1\cup\dots\cup\varpi_{j-1} )$}
    \right) ,
 \]
 hold,
 where $\mu$ and  $\nu$ are constants defined in 
 \eqref{eq:constants}.
\end{lemma}
\begin{proof}
 By \ref{key:concl:0}, \ref{key:concl:2},
    \eqref{eq:constants} and \eqref{eq:N1}, 
 \begin{align*}
   |\varphi_{j-1}|&\geq
    |\varphi_0|-|\varphi_{1}-\varphi_0|-\dots-
    |\varphi_{j-1}-\varphi_{j-2}|\\
     &\geq \min_{\overline{\D}_1}|\varphi_0|-\frac{(j-1)\varepsilon}{2N^{2}}
     \geq \nu-\frac{\varepsilon}{N}
     \geq \frac{\nu}{2}
 \end{align*}
 holds on 
  $\overline{\D}_1\setminus(\varpi_1\cup\dots\cup\varpi_{j-1} )$.
 On the other hand, 
 we have
 \begin{align*}
   |\varphi_{j-1}|&\leq
    |\varphi_0|+|\varphi_{1}-\varphi_0|+\dots
    +|\varphi_{j-1}-\varphi_{j-2}|\\
     &\leq \max_{\overline{\D}_1}|\varphi_0|+\frac{(j-1)\varepsilon}{2N^2}
     \leq \max_{\overline{\D}_1}|\varphi_0|+
           \frac{\varepsilon}{N}
              \leq \max_{\overline{\D}_1}|\varphi_0|+1\leq \mu.
  \qquad \mbox{\qed}
 \end{align*}
\renewcommand{\qed}{\relax}
\end{proof}
\begin{lemma}\label{lem:diam-omega}
 Fix an integer $j$ {\rm(}$1\le j\le 2N${\rm)}.
 If $F_0$, $F_1$, \dots $F_{j-1}$ satisfy
{ \rm \ref{key:concl:0}} and {\rm \ref{key:concl:2}} of the Key
 Lemma~\ref{lem:key}.
 Then for each $q\in\overline{\varpi}_j$, it holds that
 \[
     |F_{j-1}(q)-F_{j-1}(\zeta_j)|\leq \frac{6\mu}{N},\qquad
     |\varphi_{j-1}(q)-\varphi_{j-1}(\zeta_j)|\leq
              |\frac{6\mu+2\varepsilon}{N},
 \]
 where $\zeta_j$ is the ``base point'' of $\varpi_j$,
 see \eqref{eq:zeta} in Appendix~\ref{app:labyrinth}.
\end{lemma}
\begin{proof}
 By Lemma~\ref{lem:omega-path} in Appendix~\ref{app:labyrinth},
 there exists a path $\gamma$ in $\overline{\varpi}_j$
 joining $\zeta_j$ and $q$ 
 such that $\Length_{\C}(\gamma)\leq 6/N$.
 Since the image of $\gamma$ lies on 
 $\overline{\D}_1\setminus\bigl(\varpi_1\cup\dots\cup\varpi_{j-1}\bigr)$,
 Lemma~\ref{lem:bound-j} implies that
 we have
 \begin{align*}
  \left|F_{j-1}(q)-F_{j-1}(\zeta_j)\right|
  & \leq  
       \int_{\gamma}|\varphi_{j-1}(z)|\,|dz|
  \leq \mu\cdot\Length_{\C}(\gamma)
  \leq \frac{6\mu}{N}.
 \end{align*}
 On the other hand,
 \begin{alignat*}{2}
  |&\varphi_{j-1}(q)-\varphi_{j-1}(\zeta_j)|\\
   &\leq
    |\varphi_{j-1}(q)-\varphi_{j-2}(q)|+\dots +
     |\varphi_1(q)-\varphi_0(q)|
     +
    |\varphi_0(q)-\varphi_0(\zeta_j)|\\
   &\hphantom{\leq|\varphi_{j-1}(q)} +
    |\varphi_{j-1}(\zeta_j)-\varphi_{j-2}(\zeta_j)|+\dots +
     |\varphi_1(\zeta_j)-\varphi_0(\zeta_j)|\\
   &\leq
    \frac{2(j-1)\varepsilon}{2N^2}+|\varphi_0(q)-\varphi_0(\zeta_j)|
    \leq \frac{2\varepsilon}{N}+
         \left|
            \int_{\gamma}\varphi_0'(z)\,dz
         \right|
   \qquad&&\text{(by \ref{key:concl:2})}
  \\
   &\leq \frac{2\varepsilon}{N}+
            \int_{\gamma}|\varphi_0'(z)| |dz|
   \leq \frac{2\varepsilon}{N}+
            \mu \Length_{\C}(\gamma)
   \leq 
           \frac{2\varepsilon}{N}+\frac{6\mu}{N}
    \qquad&&\text{(by \eqref{eq:constants})}.
\qed
\end{alignat*}%
\renewcommand{\qed}{\relax}
\end{proof}
We fix $j$ ($1\leq j\leq 2N$) and assume 
$F_0$, $F_1$, \dots $F_{j-1}$ are already constructed
and satisfy \ref{key:concl:0}--\ref{key:concl:6}.
From now on, we give a recipe of construction of  
$F_j$ and $\vect{v}_j$ as an inductive procedure:

\begin{lemma}\label{lem:axis-vector}
 There exists a unit vector $\vect{u}\in\C^3$ 
 {\rm(}i.e.\ $|\vect{u}|=1${\rm)} such that
 \begin{enumerate}
  \item\label{axis:2}
       $\delta^2:=|\inner{\vect{u}}{\vect{u}}|\geq 1/{{N^{1/4}}}$.
  \item\label{axis:3}
       If 
       \begin{equation}\label{eq:zeta-j-far}
	|F_{j-1}(\zeta_j)|\geq \frac{1}{\sqrt{N}},
       \end{equation}
       it holds that
       \[
	   \left|\hinner{\frac{F_{j-1}(p)}{|F_{j-1}(p)|}}{\vect{u}}\right|
           \geq 1-\frac{c_1}{\sqrt{N}}
           \qquad (p\in\overline{\varpi}_j),
       \]
       where $c_1$ is the constant as in \eqref{eq:c123}.
 \end{enumerate}
\end{lemma}
\begin{proof}
 When $|F_{j-1}(\zeta_j)|<1/\sqrt{N}$, the unit vector 
 $\vect{u}=(0,0,1)$ satisfies the conclusions.
 (Note that the conclusion \ref{axis:3} is empty in this case.)

 Now, we assume \eqref{eq:zeta-j-far}, and set
 \begin{equation}\label{eq:v-zero}
  \vect{u}_0:= \frac{F_{j-1}(\zeta_j)}{|F_{j-1}(\zeta_j)|}.
 \end{equation}
 By Lemma~\ref{lem:diam-omega}, 
 \begin{equation}\label{eq:diam-omega}
     |F_{j-1}(p)-F_{j-1}(\zeta_j)|\leq \frac{6\mu}{N}
 \end{equation}
 holds for each $p\in \overline{\varpi}_j$.
  Then, for $p\in \overline{\varpi}_j$,
 it holds that
 \begin{alignat*}{2}
    |&F_{j-1}(p)|\geq |F_{j-1}(\zeta_j)|-|F_{j-1}(p)-F_{j-1}(\zeta_j)|
           \geq|F_{j-1}(\zeta_j)|-\frac{6\mu}{N}
         \qquad &&\text{(by \eqref{eq:diam-omega})}\\
        &\geq 
        |F_{j-1}(\zeta_j)|\left(1-\frac{6\mu}{N|F_{j-1}(\zeta_j)|}\right)
        \geq
        |F_{j-1}(\zeta_j)|\left(1-
            \frac{6\mu}{N \frac{1}{\sqrt{N}}}\right)
         \qquad &&\text{(by \eqref{eq:zeta-j-far})}\\
      &
  = |F_{j-1}(\zeta_j)| \left(1-\frac{6\mu}{\sqrt{N}}\right)
       \geq \frac{1}{2}|F_{j-1}(\zeta_j)|
      \qquad&&\text{(by \eqref{eq:N1})}.
 \end{alignat*}
 Thus, using \eqref{eq:zeta-j-far} again, we have
 \begin{equation}\label{eq:fp-lower}
    |F_{j-1}(p)|\geq \frac{1}{2}|F_{j-1}(\zeta_j)|\geq
     \frac{1}{2\sqrt{N}}
     \qquad (p \in \overline{\varpi}_j).
 \end{equation}
 Then by the relationship of the arithmetic mean and the geometric mean,
 we have
\allowdisplaybreaks{
 \begin{align*}
  &\text{%
    \begin{minipage}{0.9\textwidth}
     $\dfrac{(6\mu)^2}{N^2}\geq
       |F_{j-1}(p)-F_{j-1}(\zeta_j)|^2$
     \hfill
     (by \eqref{eq:diam-omega})
  \end{minipage}}\\
  &= |F_{j-1}(p)|^2 + |F_{j-1}(\zeta_j)|^2 - 
       2\Re\hinner{F_{j-1}(p)}{F_{j-1}(\zeta_j)}\\
  &=|F_{j-1}(p)|\,|F_{j-1}(\zeta_j)|
    \left(
     \frac{|F_{j-1}(p)|}{|F_{j-1}(\zeta_j)|}+
     \frac{|F_{j-1}(\zeta_j)|}{|F_{j-1}(p)|}-
     2\Re\hinner{
       \frac{F_{j-1}(p)}{|F_{j-1}(p)|}}{
       \frac{F_{j-1}(\zeta_j)}{|F_{j-1}(\zeta_j)|}}
     \right)\\
  &\geq
    2|F_{j-1}(p)|\,|F_{j-1}(\zeta_j)|
    \left(
     1-
     \Re\hinner{
       \frac{F_{j-1}(p)}{|F_{j-1}(p)|}
     }{
       \frac{F_{j-1}(\zeta_j)}{|F_{j-1}(\zeta_j)|}
     }
     \right)\\
  &\text{%
    \begin{minipage}{0.9\textwidth}
     $\displaystyle{\geq \frac{1}{N}
      \left(
       1-
       \Re\hinner{
        \frac{F_{j-1}(p)}{|F_{j-1}(p)|}
       }{
        \frac{F_{j-1}(\zeta_j)}{|F_{j-1}(\zeta_j)|}
       }
       \right)}$
    \hfill(by \eqref{eq:fp-lower}, \eqref{eq:zeta-j-far})
     \end{minipage}}\\
  &\geq \frac{1}{N}
    \left(
     1-
     \left|\hinner{
       \frac{F_{j-1}(p)}{|F_{j-1}(p)|}
     }{
       \frac{F_{j-1}(\zeta_j)}{|F_{j-1}(\zeta_j)|}
     }
     \right|
     \right).
 \end{align*}}%
 Hence, by \eqref{eq:N1}, we have
 \begin{equation} \label{eq:angle-est-0}
 \begin{aligned}
   \left|\hinner{
       \frac{F_{j-1}(p)}{|F_{j-1}(p)|}
     }{ \vect{u}_0
     }
     \right|&=
   \left|\hinner{
       \frac{F_{j-1}(p)}{|F_{j-1}(p)|}
     }{
       \frac{F_{j-1}(\zeta_j)}{|F_{j-1}(\zeta_j)|}
     }
     \right|\\
    &\geq 1-\frac{(6\mu)^2}{N}
    = 1 -\frac{6\mu^2}{\sqrt{N}}\frac{6}{\sqrt{N}}
    \geq 1 -\frac{6\mu^2}{\sqrt{N}}.
 \end{aligned}
 \end{equation}
\subparagraph{Case A}
We consider the case 
$|\inner{\vect{u}_0}{\vect{u}_0}|\geq 1/{{N}^{1/4}}$.
In this case, we set $\vect{u}=\vect{u}_0$.
Then the unit vector $\vect{u}$ satisfies \ref{axis:2} trivially.
Moreover, \eqref{eq:angle-est-0} implies 
the assertion \ref{axis:3} because
$c_1$ in \eqref{eq:c123} satisfies
$c_1 \geq 6\mu^2$.
\subparagraph{Case B}
We next consider the case
$|\inner{\vect{u}_0}{\vect{u}_0}|< 1/{{N}^{1/4}}$.
 In this case, set 
 \begin{equation}\label{eq:vect-v}
    \vect{u}:= 
     \frac{\tilde{\vect{u}}}{|\tilde{\vect{u}}|},\qquad
     \text{where}\quad
     \tilde{\vect{u}}:=
        \vect{u}_0+ \frac{2}{{{N}^{1/4}}}\bar{\vect{u}}_0.
 \end{equation}
 To show \ref{axis:2} and \ref{axis:3}, we set
 \begin{equation}\label{eq:delta-zero}
     \delta_0^2:=
     |\inner{\vect{u}_0}{\vect{u}_0}|\left(<\frac{1}{{{N}^{1/4}}}\right),
     \qquad (\delta_0\geq 0).
 \end{equation}
 Since $\vect{u}_0$ is a unit vector, \eqref{eq:vect-v} yields
 \begin{alignat}{2}\label{eq:u-inner-abs}
  |&\inner{\tilde{\vect{u}}}{\tilde{\vect{u}}}|
  =
  \left|
  \inner{\vect{u}_0}{\vect{u}_0}+
     \frac{4}{\sqrt{N}}
      \inner{\bar{\vect{u}}_0}{\bar{\vect{u}}_0}+
      \frac{4}{{{N}^{1/4}}}\inner{\vect{u}_0}{\bar{\vect{u}}_0}
  \right|\\
  &\geq 
  \frac{4}{{{N}^{1/4}}}|\hinner{\vect{u}_0}{\vect{u}_0}|
  -|\inner{\vect{u}_0}{\vect{u}_0}| -
  \frac{4}{\sqrt{N}}|\inner{\vect{u}_0}{\vect{u}_0}|\nonumber \\
  &\geq \frac{4}{{{N}^{1/4}}}
       -\delta^2_0\left(1+\frac{4}{\sqrt{N}}\right)
  \geq \frac{4}{{{N}^{1/4}}}-\frac{5}{3}\delta_0^2
  \geq \frac{7}{3{{N}^{1/4}}}
  \qquad &&\text{(by \eqref{eq:N1},\eqref{eq:delta-zero}).}
  \nonumber
 \end{alignat}
 On the other hand, using \eqref{eq:delta-zero} and \eqref{eq:N1} again,
 we have
 \begin{align}\label{eq:u-abs}
  |\tilde{\vect{u}}|^2&=
  |\vect{u}_0|^2 + \frac{4}{\sqrt{N}}|\bar{\vect{u}}_0|^2 +
     \frac{4}{{{N}^{1/4}}}
    \Re
      \bigl(\hinner{\vect{u}_0}{\bar{\vect{u}}_0}\bigr)\\
  &=1 + \frac{4}{\sqrt{N}}+
     \frac{4}{{{N}^{1/4}}}
     \Re \inner{\vect{u}_0}{\vect{u}_0}
  \leq 1 + \frac{4}{\sqrt{N}}+
     \frac{4}{{{N}^{1/4}}}
       |\inner{\vect{u}_0}{\vect{u}_0}|\nonumber\\
  &= 1  + \frac{4}{\sqrt{N}}+
         \frac{4\delta_0^2}{{{N}^{1/4}}}
  \leq 1  + \frac{4}{\sqrt{N}}+
         \frac{4}{\sqrt{N}}
  \leq 1+ \frac{8}{\sqrt{N}}\leq \frac{7}{3}.\nonumber
 \end{align}
 Then \ref{axis:2} holds because of \eqref{eq:u-inner-abs} and \eqref{eq:u-abs}.

 Finally, we prove \ref{axis:3}.
 Let $p\in \overline{\varpi}_j$.
 Then we have
 \begin{alignat*}{2}
    &|\inner{F_{j-1}(p)}{\vect{u}_0}|
    =
    |\inner{F_{j-1}(p)-F_{j-1}(\zeta_j)}{\vect{u}_0}
            +\inner{F_{j-1}(\zeta_j)}{\vect{u}_0}|\\
    &\leq 
    |\inner{F_{j-1}(p)-F_{j-1}(\zeta_j)}{\vect{u}_0}|
          +|\inner{F_{j-1}(\zeta_j)}{\vect{u}_0}|
    \\
    &\leq |F_{j-1}(p)-F_{j-1}(\zeta_j)|\,|\vect{u}_0| +
    \bigl|\inner{|F_{j-1}(\zeta_j)|\vect{u}_0}{\vect{u}_0}\bigr|
    \ && \text{(by \eqref{eq:schwarz}, \eqref{eq:v-zero})}\\
    &\leq \frac{6\mu}{N}+|F_{j-1}(\zeta_j)|\,|\inner{\vect{u}_0}{\vect{u}_0}|
     \qquad &&\text{(by \eqref{eq:diam-omega})}
    \\
    &\leq \frac{6\mu}{N}+\delta_0^2|F_{j-1}(\zeta_j)|
    \leq \frac{6\mu}{N}+\frac{|F_{j-1}(\zeta_j)|}{{{N}^{1/4}}}
      &&\text{(by \eqref{eq:delta-zero})}\\
    &=
      |F_{j-1}(\zeta_j)|\left(\frac{1}{{{N}^{1/4}}}+\frac{6\mu}{
            N|F_{j-1}(\zeta_j)|}\right)
\leq 
      |F_{j-1}(\zeta_j)|\left(\frac{1}{{{N}^{1/4}}}+\frac{6\mu}{\sqrt{N}}\right)
    ~
    &&\text{(by \eqref{eq:zeta-j-far})}\\
    &=|F_{j-1}(\zeta_j)|\frac{1}{{{N}^{1/4}}}
            \left(1+\frac{6\mu}{{{N}^{1/4}}}\right)
    \leq 
      |F_{j-1}(\zeta_j)|\frac{1}{{{N^{1/4}}}}
            \left(1+\frac{6\mu}{\sqrt{6}}\right)~
    &&\text{(by \eqref{eq:N1}).}
 \end{alignat*}
 Thus we have
 \begin{equation}\label{eq:f-j-u0}
    |\inner{F_{j-1}(p)}{\vect{u}_0}|
       \leq\frac{|F_{j-1}(\zeta_j)|}{{{N^{1/4}}}}(1+3\mu).
 \end{equation}
 On the other hand, since
 \begin{equation}\label{eq:sqrt-est}
     \frac{1}{\sqrt{1+x}}\geq 1-\frac{x}{2}\qquad (0\leq x\leq 2),
 \end{equation}
 we have
 \begin{alignat*}{2}
   &\left|\hinner{\frac{F_{j-1}(p)}{|F_{j-1}(p)|}}{\vect{u}}\right|
    =\frac{1}{|\tilde{\vect{u}}|}
        \left|
          \hinner{\frac{F_{j-1}(p)}{|F_{j-1}(p)|}}{\vect{u}_0}
          +\frac{2}{{{N^{1/4}}}}
          \hinner{\frac{F_{j-1}(p)}{|F_{j-1}(p)|}}{\bar{\vect{u}}_0}
        \right|~
  && \text{(by \eqref{eq:vect-v})}\\
   &\geq \frac{1}{\sqrt{1+\frac{8}{\sqrt{N}}}}
        \left|
          \hinner{\frac{F_{j-1}(p)}{|F_{j-1}(p)|}}{\vect{u}_0}
          +\frac{2}{{{N^{1/4}}}}
          \hinner{\frac{F_{j-1}(p)}{|F_{j-1}(p)|}}{\bar{\vect{u}}_0}
        \right|
    \qquad &&\text{(by \eqref{eq:u-abs})}\\
   &\geq \left(1-\frac{4}{\sqrt{N}}\right)
         \left(
          \left|
          \hinner{\frac{F_{j-1}(p)}{|F_{j-1}(p)|}}{\vect{u}_0}\right|
          -
          \frac{2}{{{N^{1/4}}}}
          \left|
          \hinner{\frac{F_{j-1}(p)}{|F_{j-1}(p)|}}{\bar{\vect{u}}_0}\right|
         \right)
    \qquad &&\text{(by \eqref{eq:sqrt-est})}\\
    &\geq \left(1-\frac{4}{\sqrt{N}}\right)
          \left[
          \left(1-\frac{(6\mu)^2}{N}\right)
          -\frac{2}{{{N^{1/4}}}}
          \left|
          \hinner{\frac{F_{j-1}(p)}{|F_{j-1}(p)|}}{\bar{\vect{u}}_0}\right|
          \right]
     \qquad &&\text{(by \eqref{eq:angle-est-0})}\\
    &\geq \left(1-\frac{4}{\sqrt{N}}\right)
          \left[
          \left(1-\frac{(6\mu)^2}{N}\right)
          -\frac{2}{{{N^{1/4}}}}
          \left|
          \inner{\frac{F_{j-1}(p)}{|F_{j-1}(p)|}}{\vect{u}_0}\right|
          \right]
     \qquad \\
    &\geq \left(1-\frac{4}{\sqrt{N}}\right)
          \left[
          \left(1-\frac{(6\mu)^2}{N}\right)
          -\frac{2|F_{j-1}(\zeta_j)|}{\sqrt{N}|F_{j-1}(p)|}(1+3\mu)
          \right]
     \qquad &&\text{(by \eqref{eq:f-j-u0})}\\
    &\geq \left(1-\frac{4}{\sqrt{N}}\right)
          \left[
          \left(1-\frac{(6\mu)^2}{N}\right)
          -\frac{4}{\sqrt{N}}(1+3\mu)
          \right]
     \qquad &&\text{(by \eqref{eq:fp-lower})}\\
    &\geq \left(1-\frac{4}{\sqrt{N}}\right)
          \left(1-\frac{1}{\sqrt{N}}\left(\frac{(6\mu)^2}{\sqrt{N}}
              +4(3\mu+1)\right)\right)
     \qquad &&\\
    &\geq \left(1-\frac{4}{\sqrt{N}}\right)
          \left(1-\frac{1}{\sqrt{N}}\left(6\mu^2
              +12\mu+4\right)\right)
     \qquad &&\text{(by \eqref{eq:N1})}\\
    &\geq   1-\frac{1}{\sqrt{N}}(6\mu^2+12\mu+8)
                 +\frac{4}{N}\left(6\mu^2+12\mu+4\right)\\
    &\geq
          1-\frac{1}{\sqrt{N}}\left(6\mu^2+12\mu+8\right)=1-\frac{c_1}{\sqrt{N}}.
    \qquad
    &&\text{(by \eqref{eq:c123})}
 \end{alignat*}
 Thus we have the conclusion.
\end{proof}

\section{The proof of the Key Lemma~\ref{lem:key}}
\label{sec:proof-key}
To continue the procedure of the iterational construction
of $F_j$, we prepare the following lemma:
\begin{lemma}\label{lem:vect-canonical}
 For a unit vector $\vect{u}\in\C^3$,
 there exists  $P\in \SO(3)$ and  $\tau\in\R$ such that 
 \begin{equation}\label{eq:P}
      e^{-\imag\tau}P\vect{u}=\pmt{0 \\ \imag \sin\theta\\ \cos\theta}
       \qquad (\imag=\sqrt{-1}).
 \end{equation}
 Here, $\theta$ is a real number such that 
 \[
   \cos 2\theta=\cos^2\theta-\sin^2\theta=|\inner{\vect{u}}{\vect{u}}|
     \qquad\left(0\leq \theta\leq \frac{\pi}{4}\right).
 \]
\end{lemma}
\begin{proof}
 Write $\vect{u}=\vect{x}+\imag \vect{y}$ ($\vect{x}$,
 $\vect{y}\in\R^3$),
 and let $\tau\in\R$ be
 \[
     \tau = 
       \begin{cases}
	\dfrac{1}{2}\arctan \dfrac{2\hinner{\vect{x}}{\vect{y}}}{%
              |\vect{x}|^2-|\vect{y}|^2}
	\qquad
	&\text{%
	(when $|\vect{x}|\neq |\vect{y}|$)}\\[10pt]
        \dfrac{\pi}{4}\qquad
	&\text{%
	(when $|\vect{x}|=|\vect{y}|$)}.
       \end{cases}
 \]
 Then $\tilde{\vect{u}}:=e^{-\imag\tau}\vect{u}$
 satisfies
 $\hinner{\Re\tilde{\vect{u}}}{\Im\tilde{\vect{u}}}=0$.
 Moreover, 
 replacing $\tau$ with $\tau+\frac{\pi}{2}$ if necessary,
 we may assume
 \begin{equation}\label{eq:re-im-2}
     |\Re\tilde{\vect{u}}|\geq |\Im\tilde{\vect{u}}|
 \end{equation}
 without loss of generality.
 In particular, since 
 $|\tilde{\vect{u}}|=1$,
 it holds that $|\Re\tilde{\vect{u}}|>0$.
 Hence there exists a matrix $P_1\in\SO(3)$ such that
 \[
     P_1(\Re\tilde{\vect{u}})=\pmt{0\\ 0\\ t}\qquad (t>0).
 \]
 Since $P_1$ is a real matrix, $\Im(P_1\tilde{\vect{u}})$
 is orthogonal to $\Re(P_1\tilde{\vect{u}})$.
 Hence, we have
 \[
     P_1(\tilde{\vect{u}})=\pmt{\imag u_1\\ \imag u_2\\ t}
     \qquad
      (u_1,u_2,t\in\R, t>0, (u_1)^2+(u_2)^2+t^2=1).
 \]
 Moreover, 
 $t \geq \sqrt{(u_1)^2+(u_2)^2}$  holds because \eqref{eq:re-im-2}.
 Next, choose a real number $s$ such that
 \[
    \begin{pmatrix}
     \cos s & -\sin s \\
     \sin s &\hphantom{-}\cos s
    \end{pmatrix}
    \begin{pmatrix}
     u_1 \\ u_2
    \end{pmatrix} = \begin{pmatrix}
		      0 \\ u
		    \end{pmatrix}
    \qquad u=\sqrt{(u_1)^2+(u_2)^2}\geq 0,
 \]
 and set
 \[
    P:=
    \begin{pmatrix}
     \cos s & -\sin s & 0\\
     \sin s &\hphantom{-}\cos s& 0 \\
      0 & 0 & 1
     \end{pmatrix} \cdot P_1\in \SO(3).
 \]
 Then 
 \[
    e^{-\imag \tau} P\vect{u}
     = P\tilde{\vect{u}} = 
     \begin{pmatrix}0\\\imag u \\ t\end{pmatrix}
    \qquad (u,t\in\R, t\geq u\geq 0, u^2 + t^2=1).
 \]
 Hence there exists $\theta\in[0,\frac{\pi}{4}]$
 such that $u=\sin\theta$, $t=\cos\theta$.
 In particular, 
 \[
    |\inner{\vect{u}}{\vect{u}}|
    = t^2 - u^2
    = \cos^2\theta-\sin^2\theta=\cos 2\theta
 \]
 holds and thus we have the conclusion.
\end{proof}
We set
\begin{equation}\label{eq:matA}
   A:= \begin{pmatrix}
     \sqrt{\cos2\theta} & 0 & 0 \\
     0 & \cos\theta & -\imag\sin\theta \\
     0 & \imag\sin\theta & \cos\theta
    \end{pmatrix},
\end{equation}
where  $\theta\in[0,\frac{\pi}{4}]$.
Then $A$ is non-singular if and only if 
$\theta\neq \frac{\pi}{4}$.
In this case,
\begin{equation}\label{eq:delta-A}
  A \in\delta\cdot \O(3,\C)\qquad (\delta:=\sqrt{\cos 2\theta}),
\end{equation}
and 
\begin{equation}\label{eq:matAinv}
      A^{-1}=\frac{1}{\delta^2}
           \begin{pmatrix}
	    \sqrt{\cos2\theta} & 0 & 0 \\
	    0 & \cos\theta & \imag\sin\theta \\
	    0 & -\imag\sin\theta & \cos\theta
	   \end{pmatrix}.
\end{equation}
\begin{lemma}\label{lem:o3-matrix-norm}
 Let $\theta\in[0,\frac{\pi}{4})$ be a real number.
 Then the matrix $A$ in \eqref{eq:matA} satisfies
 \[
    \Vert{}A\Vert{} = \cos\theta+\sin\theta \leq \sqrt{2},\qquad
    \Vert{}A^{-1}\Vert{} = \frac{\cos\theta+\sin\theta}{\cos 2\theta}
          \leq \frac{\sqrt{2}}{\delta^2},\qquad
           (\delta=\sqrt{\cos 2\theta}).
 \]
\end{lemma}
\begin{proof}
 Since the eigenvalues of the matrix
 \[
     A^*A
      = 
      \begin{pmatrix}
       \cos 2\theta & 0 & 0 \\
       0 & 1 & -\imag\sin 2\theta \\
       0 & \imag\sin 2\theta & 1
      \end{pmatrix}
 \]
 are
 $(\cos\theta-\sin\theta)^2$,
 $\cos 2\theta$, and $(\cos\theta+\sin\theta)^2$,
 \eqref{eq:matrix-norm-1} implies that
 $\Vert{}A\Vert{} = \cos\theta+\sin\theta
          = \sqrt{2}\sin\left(\theta+\frac{\pi}{4}\right)
          \leq \sqrt{2}$.
 On the other hand, 
 the eigenvalues of $(A^{-1})^*A^{-1}$ are
 \[
    \frac{(\cos\theta-\sin\theta)^2}{\cos^2 2\theta},\qquad
    \frac{1}{\cos 2\theta},\quad\text{and}\quad
    \frac{(\cos\theta+\sin\theta)^2}{\cos^2 2\theta}.
 \]
 Hence
 \[
    \Vert{}A^{-1}\Vert{} = \frac{\cos\theta+\sin\theta}{\cos 2\theta}
          \leq \frac{\sqrt{2}}{\delta^2},
 \]
 which is the conclusion.
\end{proof}
We return to the construction of $F_j$:
Take $\vect{u}$ as in Lemma~\ref{lem:axis-vector},
and take $P\in\SO(3)$ and $\tau\in\R$  as in Lemma~\ref{lem:vect-canonical}.
where $\theta\in [0,\frac{\pi}{4}]$ is given by 
$\cos 2 \theta=|\inner{\vect{u}}{\vect{u}}|$.
Observe that by \ref{axis:2} of Lemma~\ref{lem:axis-vector} we have
\begin{equation}\label{eq:delta-2}
   \delta:=\sqrt{\cos 2 \theta} \geq \frac{1}{{{N^{1/4}}}}
\end{equation}
and therefore $\theta\in \left[0,\frac{\pi}{4}\right)$.
We set 
\begin{equation}\label{eq:f}
   F:= e^{\imag\tau} P F_{j-1},\qquad
   \varphi := \varphi_F=\frac{dF}{dz}=e^{\imag\tau}P\varphi_{j-1}.
\end{equation}
Since $P\in\SO(3)\subset \O(3,\C)$, $F$ is a holomorphic
null immersion.
On the other hand, since $P\in\SO(3)\subset\U(3)$, 
$F$ is congruent to $F_{j-1}$ in $\C^3$.
In particular, 
\begin{equation}\label{eq:phi-abs}
   |\varphi|=|\varphi_{j-1}|,\qquad
   |\varphi(q)-\varphi(p)|=|\varphi_{j-1}(q)-\varphi_{j-1}(p)|
\end{equation}
hold for $p,q\in\overline{\D}_1$.

Taking into account \eqref{eq:delta-2}, we consider the matrix
\begin{equation}\label{eq:matrix-A}
    A=(\vect{a}^{(1)},\vect{a}^{(2)},\vect{a}^{(3)})
     := \begin{pmatrix}
	\delta & 0 & 0 \\
	0  & \cos\theta & -\imag\sin\theta\\
	0  & \imag\sin\theta & \cos\theta
     \end{pmatrix}
       \in \delta\cdot \O(3,\C).
\end{equation}
In particular, by \eqref{eq:P}, it holds that
\begin{equation}\label{eq:a-3}
   \vect{a}^{(3)} = \overline{e^{-\imag\tau}P\vect{u}}=
        e^{\imag\tau}P\bar{\vect{u}}.
\end{equation}
By Lemma~\ref{lem:o3-matrix-norm}  and 
\eqref{eq:delta-2}, it holds that
\begin{equation}\label{eq:a-norm}
   \Vert{}A\Vert{}\leq \sqrt{2},\qquad 
   \Vert{}A^{-1}\Vert{}\leq \frac{\sqrt{2}}{\delta^2}
                                      \leq \sqrt{2}{{N^{1/4}}}.
\end{equation}
Using the matrix $A$ in \eqref{eq:matrix-A}, we set 
\begin{equation}\label{eq:cal-f}
 \begin{aligned}
   E&=\trans{(E^{(1)},E^{(2)},E^{(3)})} 
    := A^{-1}F = e^{\imag\tau}A^{-1}PF_{j-1},\\
   \psi&:= \frac{dE}{dz}=A^{-1}\varphi
         =e^{\imag \tau} A^{-1}P\varphi_{j-1}.
 \end{aligned}
\end{equation}
Since  $A\in\delta\cdot\O(3,\C)$,
$E$ is a holomorphic null immersion although it is not necessarily congruent 
to $F_{j-1}$.
Moreover, by
\eqref{eq:cal-f}, 
\eqref{eq:norm-ineq},
\eqref{eq:a-norm} and \eqref{eq:phi-abs}, we have
\begin{align}\label{eq:psi-below}
   |\psi|&= |A^{-1}\varphi| \geq \frac{1}{\Vert{}A\Vert{}}|\varphi|\geq 
          \frac{|\varphi|}{\sqrt{2}}=
          \frac{|\varphi_{j-1}|}{\sqrt{2}}\\
    \label{eq:psi-above}
   |\psi(q)-\psi(p)|&=
   |A^{-1}\bigl(\varphi(q)-\varphi(p)\bigr)|
   \leq \Vert{}A^{-1}\Vert{}\,|\varphi(q)-\varphi(p)|\\
    &\leq \sqrt{2}{{N^{1/4}}}\,|\varphi_{j-1}(q)-\varphi_{j-1}(p)|.
 \nonumber
\end{align}

\begin{lemma}\label{lem:rotation}
 Let $G=G_E\colon{}\overline{\D}_1\to S^2$ 
 be the Gauss map of $E$ as in \eqref{eq:cal-f}
 {\rm(}cf.\ \eqref{eq:gauss}{\rm)}.
 Then there exists a real matrix $Q$
 \begin{equation}\label{eq:matrix-Q}
    Q = \begin{pmatrix}
	 1 & 0 & 0 \\
	 0 & \cos\Theta & -\sin\Theta \\
	 0 & \sin\Theta & \hphantom{-}\cos\Theta
	\end{pmatrix}\in \SO(3),
    \qquad |\Theta|\leq \frac{4}{\sqrt{N}}
 \end{equation}
 such that 
 \begin{equation}\label{eq:gauss-dist}
      \dist_{S^2}\bigl(Q G(p),\pm\vect{e}_3\bigr)\geq \frac{1}{\sqrt{N}}
      \qquad (\vect{e}_3=(0,0,1))
 \end{equation}
 holds for each point $p\in\overline{\varpi}_j$,
 where
 $\dist_{S^2}$ is the canonical distance function 
 of the unit sphere $S^2$.
 In particular, 
 the matrix $Q$ 
 commutes with $A^{-1}$ as in \eqref{eq:matrix-A},
 that is, 
 \begin{equation}\label{eq:commute}
AQ^{-1}=Q^{-1}A,
 \end{equation}
 and 
 \begin{equation}\label{eq:norm-q}
   \Vert{}Q^{-1}-\id\Vert{}\leq |\Theta|\leq \frac{4}{\sqrt{N}}
 \end{equation}
 holds.
\end{lemma}
\begin{proof}
 By \eqref{eq:matrix-A} and \eqref{eq:matrix-Q},
 \eqref{eq:commute} is trivial.
 Moreover, since
 \[
   Q^{-1}-\id 
   = \begin{pmatrix}
      0 & 0 & 0 \\
      0 & \cos\Theta -1 & \sin\Theta \\
      0 & -\sin\Theta   & \cos\Theta -1
     \end{pmatrix}
   = -2\sin\frac{\Theta}{2}
     \begin{pmatrix}
      0 & 0& 0 \\
      0 & \sin\frac{\Theta}{2} & -\cos\frac{\Theta}{2}\\
      0 & \cos\frac{\Theta}{2} & \hphantom{-}\sin\frac{\Theta}{2}
     \end{pmatrix},
 \]
 the maximum eigenvalue of 
 $(Q^{-1}-\id)^*(Q^{-1}-\id)$ is
 $\left(2\sin\frac{\Theta}{2}\right)^2$.
 Hence by \eqref{eq:matrix-norm-1},
 \[
   \Vert{}Q^{-1}-\id\Vert{}=2\left|\sin\frac{\Theta}{2}\right|\leq |\Theta|.
 \]
 holds, and thus we have \eqref{eq:norm-q}.

 So it is sufficient to show \eqref{eq:gauss-dist}
 for suitable $Q$. The Euclidean distance between $G(p)$ and $G(\zeta_j)$ in $\R^3$
 can be estimated as 
\allowdisplaybreaks{%
 \begin{align*}
  &
  \text{\begin{minipage}{\textwidth}$\displaystyle{
	 |G(p)-G(\zeta_j)|=
	 \left|\frac{\psi(p)\times\overline{\psi(p)}}{|\psi(p)|^2}-
        \frac{\psi(\zeta_j)\times\overline{\psi(\zeta_j)}}{|\psi(\zeta_j)|^2}
        \right|}$\hfill
  (by \eqref{eq:gauss})
  \end{minipage}}\\
  &=
  \frac{1}{|\psi(p)|^2|\psi(\zeta_j)|^2}
  \left|
     \left(\psi(p)\times\overline{\psi(p)}\right)|\psi(\zeta_j)|^2 - 
     \left(\psi(\zeta_j)\times\overline{\psi(\zeta_j)}\right)|\psi(p)|^2 
  \right|\\
  &=
  \frac{1}{|\psi(p)|^2|\psi(\zeta_j)|^2}
  \left|
     \left(\psi(p)\times\overline{\psi(p)}\right)|\psi(\zeta_j)|^2 - 
     \left(\psi(p)\times\overline{\psi(p)}\right)|\psi(p)|^2 \right.\\
  &\left.\hspace{0.3\textwidth}
    +
     \left(\psi(p)\times\overline{\psi(p)}\right)|\psi(p)|^2 -
     \left(\psi(\zeta_j)\times\overline{\psi(\zeta_j)}\right)|\psi(p)|^2 
  \right|\\
  &= 
  \frac{1}{|\psi(p)|^2|\psi(\zeta_j)|^2}
  \left|
     \left(\psi(p)\times\overline{\psi(p)}\right)
    \left(|\psi(\zeta_j)|^2 - |\psi(p)|^2\right)\right.\\
  &\left.\hspace{0.3\textwidth}
     +
     |\psi(p)|^2 \left(
     \psi(p)\times\overline{\psi(p)}-
     \psi(\zeta_j)\times\overline{\psi(\zeta_j)}\right)
  \right|\\
  &
  \text{\begin{minipage}{\textwidth}$\displaystyle{
	 \leq 
	 \frac{|\psi(p)|^2
	 \left(
	 \left||\psi(\zeta_j)|^2 - |\psi(p)|^2\right|
	 +\left|
	 \psi(p)\times\overline{\psi(p)}-
	 \psi(\zeta_j)\times\overline{\psi(\zeta_j)}\right|\right)
	 }{|\psi(p)|^2|\psi(\zeta_j)|^2}}$
	 \hfill
        ($|\varphi\times\bar\varphi|=|\varphi|^2$)
	\end{minipage}}\\
  &=
  \frac{1}{|\psi(\zeta_j)|^2}
  \left(
     \left||\psi(\zeta_j)|^2 - |\psi(p)|^2\right|
     +\left|
     \psi(p)\times\overline{\psi(p)}-
     \psi(\zeta_j)\times\overline{\psi(\zeta_j)}\right|\right)\\
  &\leq 
  \frac{1}{|\psi(\zeta_j)|^2}
  \biggl(
  \bigl|
     |\psi(\zeta_j)| - |\psi(p)|
  \bigr|
  \bigl(
     |\psi(\zeta_j)| + |\psi(p)|
  \bigr)\\
  &\hspace{0.2\textwidth}+
    \left|
     \psi(p)\times\overline{\psi(p)}-\psi(p)\times\overline{\psi(\zeta_j)}+
     \psi(p)\times\overline{\psi(\zeta_j)}-
     \psi(\zeta_j)\times\overline{\psi(\zeta_j)}\right|\biggr)\\
  &\leq 
  \frac{1}{|\psi(\zeta_j)|^2}
  \biggl(
  \bigl|
     \psi(\zeta_j) - \psi(p)
  \bigr|
  \bigl(
     |\psi(\zeta_j)| + |\psi(p)|
  \bigr)\\
  &\hspace{0.3\textwidth}+
    \left|
     \psi(p)\times\left(\overline{\psi(p)-\psi(\zeta_j)}\right)+
     \biggl(\psi(p)-\psi(\zeta_j)\biggr)\times\overline{\psi(\zeta_j)}
     \right|\biggr)\\
   & \leq      
      \frac{2}{|\psi(\zeta_j)|^2}
        \bigl|\psi(p)-\psi(\zeta_j)\bigr|
        \bigl(|\psi(p)-\psi(\zeta_j)| + 2|\psi(\zeta_j)|\bigr)\\
    &\leq      
      \frac{2}{|\psi(\zeta_j)|}
        \bigl|\psi(p)-\psi(\zeta_j)\bigr|
        \left(2+\frac{|\psi(p)-\psi(\zeta_j)|}{|\psi(\zeta_j)|}\right)\\
    &
      \text{\begin{minipage}{\textwidth}$\displaystyle{
      \leq 
      \frac{2\sqrt{2}}{|\varphi_{j-1}(\zeta_j)|}
        \bigl|\psi(p)-\psi(\zeta_j)\bigr|
        \left(2+\frac{\sqrt{2}
         |\psi(p)-\psi(\zeta_j)|}{|\varphi_{j-1}(\zeta_j)|}\right)}
       $\hfill	     
	     (by \eqref{eq:psi-below})
      \end{minipage}}\\
    &
      \text{\begin{minipage}{\textwidth}$\displaystyle{
	     \leq 
      \frac{4{{N^{1/4}}}
      \bigl|\varphi_{j-1}(p)-\varphi_{j-1}(\zeta_j)\bigr|}{%
      |\varphi_{j-1}(\zeta_j)|}
        \left(2+\frac{2{{N^{1/4}}}
         |\varphi_{j-1}(p)-\varphi_{j-1}(\zeta_j)|}{%
            |\varphi_{j-1}(\zeta_j)|}\right)   
	     }$
	     \hfill
           (by \eqref{eq:psi-above})
      \end{minipage}}\\
    &\text{\begin{minipage}{\textwidth}$\displaystyle{
      \leq 
      \frac{8{{N^{1/4}}}
      \left|\varphi_{j-1}(p)-\varphi_{j-1}(\zeta_j)\right|}{\nu}
        \left(2+\frac{4{{N^{1/4}}}}{\nu}
         |\varphi_{j-1}(p)-\varphi_{j-1}(\zeta_j)|\right)}$
      \hfill     
	     (Lemma~\ref{lem:bound-j})
     \end{minipage}}\\
    &\text{\begin{minipage}{\textwidth}$\displaystyle{
      \leq \frac{8{{N^{1/4}}}}{\nu}\frac{6\mu+2\varepsilon}{N}
       \left(2+\frac{4{{N^{1/4}}}}{\nu}\frac{6\mu+2\varepsilon}{N}\right)}
     $\hfill      
     (Lemma~\ref{lem:diam-omega})
     \end{minipage}}\\
    &\leq \frac{1}{\sqrt{N}} \frac{1}{{{N^{1/4}}}}
      \left(
	 \frac{16(3\mu+\varepsilon)}{\nu}
         \left(2+\frac{4(6\mu+\varepsilon)}{N^{3/4}\nu}\right)\right)\\
    &\text{\begin{minipage}{\textwidth}$\displaystyle{
      \leq \frac{1}{\sqrt{N}} \frac{1}{{{N^{1/4}}}}
      \left(
	 \frac{16(3\mu+\varepsilon)}{\nu}
         \left(2+\frac{4(6\mu+\varepsilon)}{\sqrt{6}^3\nu}\right)\right)}$
      \hfill(by \eqref{eq:N1})
      \end{minipage}}\\
    &\text{\begin{minipage}{\textwidth}$\displaystyle{
      \leq \frac{1}{\sqrt{N}} \frac{1}{{{N^{1/4}}}}
      \left(
	 \frac{16(3\mu+\varepsilon)}{\nu}
         \left(2+\frac{6\mu+\varepsilon}{3\nu}\right)\right)
     \leq \frac{1}{2\sqrt{N}}}$
      \hfill(by \eqref{eq:N1}).
	   \end{minipage}}
 \end{align*}}%
 Then we have
 \begin{align}\label{eq:diam-gauss}
   \dist_{S^2}\bigl(G(p),G(\zeta_j)\bigr)
   &= 2\arcsin\left( \frac{1}{2} | G(p)-G(\zeta_j)|  \right)\\
   &\leq \frac{\pi}{2}|G(p)-G(\zeta_j)|
   \leq 2|G(p)-G(\zeta_j)|\leq \frac{1}{\sqrt{N}}.\nonumber
 \end{align}
 Here we used the inequality 
 $\arcsin x \leq \pi x/2$ ($0\leq x\leq 1$).
 In particular, $G(\varpi_j)$ is contained 
 in the geodesic disc in $S^2$ centered at  $G(\zeta_j)$
 with radius $1/\sqrt{N}$.

\subparagraph{Case 1:}
 Assume both 
 $\dist_{S^2}\bigl(G(\zeta_j),\vect{e}_3\bigr)\geq 2/\sqrt{N}$
 and
 $\dist_{S^2}\bigl(G(\zeta_j),-\vect{e}_3\bigr)\geq 2/\sqrt{N}$
 hold.
 Then for each $p\in\overline{\varpi}_j$, \eqref{eq:diam-gauss}
 implies that
 \begin{align*}
  \dist_{S^2}\bigl(G(p),\vect{e}_3\bigr)&\geq
  \dist_{S^2}\bigl(G(\zeta_j),\vect{e}_3\bigr)-
  \dist_{S^2}\bigl(G(p),G(\zeta_j)\bigr)\\
  &\geq
  \frac{2}{\sqrt{N}}-
  \dist_{S^2}\bigl(G(p),G(\zeta_j)\bigr)
  \geq 
  \frac{2}{\sqrt{N}}-\frac{1}{\sqrt{N}}=\frac{1}{\sqrt{N}}.
 \end{align*}
 Similarly $\dist_{S^2}\bigl(G(p),-\vect{e}_3\bigr)\geq 1/\sqrt{N}$
 holds.
 Then we have the conclusion \eqref{eq:gauss-dist} for $Q=\id$ and
 $\Theta=0$.

\subparagraph{Case 2:}
 Assume 
\begin{equation}\label{eq:gauss-near}
\dist_{S^2}\bigl(G(\zeta_j),\vect{e}_3\bigr)<\frac{2}{\sqrt{N}}.
\end{equation}
 In this case, take the matrix $Q$ as in \eqref{eq:matrix-Q} with
 \begin{equation}\label{eq:vartheta}
     \Theta:=\frac{4}{\sqrt{N}}.
 \end{equation}
 Then
 \begin{align*}
  \dist&{}_{S^2}\bigl(QG(p),\vect{e}_3\bigr)\\
  &
  \geq 
  \dist_{S^2}\bigl(Q\vect{e}_3,\vect{e}_3\bigr)
  -\dist_{S^2}\bigl(QG(p),QG(\zeta_j)\bigr)
  -\dist_{S^2}\bigl(QG(\zeta_j),Q\vect{e}_3\bigr)\\
  &\text{\begin{minipage}{0.95\textwidth}$\displaystyle{
  =
   \frac{4}{\sqrt{N}}
  -\dist_{S^2}\bigl(QG(p),QG(\zeta_j)\bigr)
  -\dist_{S^2}\bigl(QG(\zeta_j),Q\vect{e}_3\bigr)
  }$\hfill	  
  (by \eqref{eq:vartheta})
  \end{minipage}}\\
  &\text{\begin{minipage}{0.95\textwidth}$\displaystyle{
  = \frac{4}{\sqrt{N}}
  -\dist_{S^2}\bigl(G(p),G(\zeta_j)\bigr)
  -\dist_{S^2}\bigl(G(\zeta_j),\vect{e}_3\bigr)
  }$\hfill	  
  ($Q\in\SO(3)$)
  \end{minipage}}\\
  &\text{\begin{minipage}{0.95\textwidth}$\displaystyle{
  \geq \frac{4}{\sqrt{N}}
     -\frac{1}{\sqrt{N}}
  -\dist_{S^2}\bigl(G(\zeta_j),\vect{e}_3\bigr)
  >\frac{1}{\sqrt{N}}
  }$\hfill	  
  (by \eqref{eq:diam-gauss}, \eqref{eq:gauss-near}).
  \end{minipage}}
\intertext{%
  On the other hand, 
}
  \dist&{}_{S^2}\bigl(QG(p),\vect{e}_3\bigr)\\
  &\leq \dist_{S^2}\bigl(QG(p),QG(\zeta_j)\bigr)+
        \dist_{S^2}\bigl(QG(\zeta_j),Q\vect{e}_3\bigr)+
        \dist_{S^2}\bigl(Q\vect{e}_3,\vect{e}_3\bigr)\\
  &\text{\begin{minipage}{0.95\textwidth}$\displaystyle{
   =\dist_{S^2}\bigl(QG(p),QG(\zeta_j)\bigr)+
        \dist_{S^2}\bigl(QG(\zeta_j),Q\vect{e}_3\bigr)+
        \frac{4}{\sqrt{N}}
  }$\hfill	  
  (by \eqref{eq:vartheta})
  \end{minipage}}
  \\
  &\text{\begin{minipage}{0.95\textwidth}$\displaystyle{
    = \dist_{S^2}\bigl(G(p),G(\zeta_j)\bigr)+
        \dist_{S^2}\bigl(G(\zeta_j),\vect{e}_3\bigr)+
        \frac{4}{\sqrt{N}}
  }$\hfill	  
  ($Q\in\SO(3)$)
  \end{minipage}}\\
  &\text{\begin{minipage}{0.95\textwidth}$\displaystyle{
   \leq \frac{1}{\sqrt{N}}+
        \dist_{S^2}\bigl(G(\zeta_j),\vect{e}_3\bigr)+
        \frac{4}{\sqrt{N}}
  }$\hfill	  
     (by \eqref{eq:diam-gauss})
  \end{minipage}}\\
  &\text{\begin{minipage}{0.95\textwidth}$\displaystyle{
  < \frac{1}{\sqrt{N}}+
        \frac{2}{\sqrt{N}}+
        \frac{4}{\sqrt{N}}=\frac{7}{\sqrt{N}}
  }$\hfill	  
    (by \eqref{eq:gauss-near})
  \end{minipage}}
 \end{align*}
 and then, 
 \begin{align*}
  \dist_{S^2}\bigl(QG(p),-\vect{e}_3\bigr)
       = \pi - \dist_{S^2}\bigl(QG(p),\vect{e}_3\bigr)\geq
              3 - \frac{7}{6}
       \geq \frac{1}{\sqrt{N}}
 \end{align*}
 because of \eqref{eq:N1}.
 Thus, we have the conclusion \eqref{eq:gauss-dist}.
\subparagraph{Case 3:}
 If $\dist_{S^2}\bigl(G(\zeta_j),-\vect{e}_3\bigr)<2/\sqrt{N}$
 holds, then we have the conclusion by the same way as in the 
 previous case.
\end{proof}
Using $P\in\SO(3)$, $\tau\in\R$ in \eqref{eq:P},
$A\in\delta\cdot\O(3,\C)$ in \eqref{eq:matrix-A} and $Q\in\SO(3)$
in \eqref{eq:matrix-Q}, we define
\begin{equation}\label{eq:tildeE}
  \widetilde E:= QE=  B^{-1}F_{j-1},\qquad
  \tilde\psi:= \frac{d\widetilde E}{dz}=Q\psi,
\end{equation}
where 
\begin{equation}\label{eq:matrix-B}
  B=
   \bigl(
     \vect{b}^{(1)},
     \vect{b}^{(2)},
     \vect{b}^{(3)}
   \bigr):=
   \bigl(e^{\imag\tau}QA^{-1}P\bigr)^{-1}
   \in (e^{-\imag\tau}\delta)\cdot \O(3,\C).
\end{equation}
Then $\widetilde E$ is a holomorphic null immersion
which is congruent to $E$ in \eqref{eq:cal-f}.
Denote by $(g,\eta)$ the Weierstrass data (cf.\ \eqref{eq:weierstrass})
of $\widetilde E$:
\begin{equation}\label{eq:weierstrass-f}
 \tilde\psi = \frac{1}{2}\bigl(1-g^2,\imag(1+g^2), 2 g\bigr)\eta,\qquad
  |\tilde\psi|^2 = \frac{1}{2}\bigl(1+|g|^2\bigr)^2|\eta|^2.
\end{equation}
Then we have
\begin{lemma}\label{lem:gauss-bound}
 The meromorphic function $g$ as in  \eqref{eq:weierstrass-f} 
 satisfies
 \[
    \frac{1}{2\sqrt{N}}\leq |g| \leq 2\sqrt{N}\qquad\text{and}\qquad
    \frac{|g|}{1+|g|^2}\geq \frac{2\sqrt{N}}{1+4N}
    \qquad (\text{on $\overline{\varpi}_j$}).
 \]
\end{lemma}
\begin{proof}
 The Gauss map $\widetilde G$ of  $\widetilde E$ 
 is obtained by 
 \[
    \widetilde G =Q G= 
    \frac{1}{1+|g|^2}
    \pmt{  2\Re g\\ 2\Im g\\|g|^2-1}.
 \]
 Here, since $\widetilde G=QG$ satisfies \eqref{eq:gauss-dist}
 on  $\overline{\varpi}_j$, 
 it holds that
 \begin{align}
     \dist_{S^2}(\widetilde G,\vect{e}_3)&= 
     \arccos\left(
                \widetilde G\cdot\vect{e}_3
              \right)=
     \arccos\left(\frac{|g|^2-1}{|g|^2+1}\right)\geq
     \frac{1}{\sqrt{N}},\label{eq:gauss-est-1}\\
     \dist_{S^2}(\widetilde G,-\vect{e}_3)&= 
     \arccos\left(
                \widetilde G\cdot(-\vect{e}_3)
              \right)=
     \arccos\left(\frac{1-|g|^2}{|g|^2+1}\right)\geq
      \frac{1}{\sqrt{N}}\label{eq:gauss-est-2} 
 \end{align}
 on $\overline{\varpi}_j$,
 where ``$\cdot$'' denotes the canonical inner product of $\R^3$.
 Since \eqref{eq:gauss-est-1} implies
 \[
     \frac{|g|^2-1}{|g|^2+1}\leq \cos\frac{1}{\sqrt{N}},
 \]
 we have
 \[
     |g|^2 \leq
       \frac{1+\cos\frac{1}{\sqrt{N}}}{1-\cos\frac{1}{\sqrt{N}}}
       = \cot^2\frac{1}{2\sqrt{N}}\leq (2\sqrt{N})^2.
 \]
 Similarly, by \eqref{eq:gauss-est-2},
 we have
 \[
     |g|^2 \geq \tan^2 \frac{1}{2\sqrt{N}}\geq \left(\frac{1}{2\sqrt{N}}\right)^2.
 \]
 Thus, we have the first inequality of the conclusion.
 The second inequality is obtained  immediately by the first inequality.
\end{proof}
We set 
\begin{equation}\label{eq:v-j}
 \vect{v}_j := \overline{\vect{b}^{(3)}},
\end{equation}
where $\vect{b}^{(3)}$ is the third column of the matrix $B$
as in \eqref{eq:matrix-B}.
\begin{lemma}\label{lem:c-45}
 The vector $\vect{v}_j$ in \eqref{eq:v-j} is a unit vector 
 satisfying 
 $|\inner{\vect{v}_j}{\vect{v}_j}|\geq 1/{{N^{1/4}}}$.
 Moreover, 
 when \eqref{eq:zeta-j-far} holds, that is, 
 $|F_{j-1}(\zeta_j)|\geq {1/\sqrt{N}}$, 
 it holds that
 \[
    \left|\hinner{
                  \frac{F_{j-1}(p)}{|F_{j-1}(p)|}
                }{\vect{v}_j}\right|
                \geq 1-\frac{C_2}{\sqrt{N}}
           \qquad \text{%
                   for $p\in\overline{\varpi}_j$},
 \]
 where $C_2$ is the constant in \eqref{eq:C1C2}.
\end{lemma}
\begin{proof}
 Let $\vect{e}_3=(0,0,1)$.
 Since the matrix $A$ and $Q^{-1}$ commute (cf.\ \eqref{eq:commute}), 
 the third column of the matrix $B$ is obtained as 
 \begin{alignat*}{2}
  \vect{b}^{(3)} &= B\vect{e}_3
                  =e^{-\imag\tau}P^{-1}AQ^{-1}\vect{e}_3
                 =e^{-\imag\tau}P^{-1}Q^{-1}A\vect{e}_3
                 \qquad&&\text{(by \eqref{eq:matrix-B}, \eqref{eq:commute})}\\
                 &=e^{-\imag\tau}P^{-1}Q^{-1}\vect{a}^{(3)}
                 =e^{-\imag\tau}P^{-1}Q^{-1}(e^{\imag\tau}P\bar{\vect{u}})
                  \qquad&&\text{(by \eqref{eq:matrix-A},\eqref{eq:a-3})}\\
                 &=P^{-1}Q^{-1}{P}\bar{\vect{u}}.
 \end{alignat*}
 Taking into account that $P$ and $Q$ are real matrices, \eqref{eq:v-j}
 implies that
 $\vect{v}_j = P^{-1}Q^{-1}P\vect{u}$.
 Then by Lemma~\ref{lem:axis-vector}, we have
 $|\vect{v}_j|=1$, $|\inner{\vect{v}_j}{\vect{v}_j}|\geq 1/{{N^{1/4}}}$,
 because $P$, $Q\in\SO(3)$.
 Moreover, when $|F_{j-1}(\zeta_j)|\geq 1/\sqrt{N}$ (i.e.\ 
 \eqref{eq:zeta-j-far} holds),
\allowdisplaybreaks{
 \begin{alignat*}{2}
  &\left|
   \hinner{\frac{F_{j-1}(p)}{|F_{j-1}(p)|}}{\vect{v}_j}
  \right|
   =
  \left|
   \hinner{\frac{F_{j-1}(p)}{|F_{j-1}(p)|}}{P^{-1}Q^{-1}P\vect{u}}
  \right|\\
  &=
  \left|
   \hinner{\frac{F_{j-1}(p)}{|F_{j-1}(p)|}}{%
     \vect{u}+P^{-1}(Q^{-1}-\id)P\vect{u}}
  \right|\\
   &\geq
  \left|
   \hinner{\frac{F_{j-1}(p)}{|F_{j-1}(p)|}}{%
     \vect{u}}\right|-
  \left|
   \hinner{\frac{F_{j-1}(p)}{|F_{j-1}(p)|}}{%
     P^{-1}(Q^{-1}-\id)P\vect{u}}
  \right|\\
  &\geq 
  \left|
   \hinner{\frac{F_{j-1}(p)}{|F_{j-1}(p)|}}{%
     \vect{u}}\right|-
  \left|
   \frac{F_{j-1}(p)}{|F_{j-1}(p)|}\right|\cdot
   \Vert{}P^{-1}(Q^{-1}-\id)P\Vert{}\,|\vect{u}|
  \quad
  &&\text{(by \eqref{eq:norm-ineq})}
  \\
  &=
  \left|
   \hinner{\frac{F_{j-1}(p)}{|F_{j-1}(p)|}}{%
     \vect{u}}\right|-
   \Vert{}P^{-1}(Q^{-1}-\id)P\Vert{}\\
  &=
  \left|
   \hinner{\frac{F_{j-1}(p)}{|F_{j-1}(p)|}}{%
     \vect{u}}\right|-
   \Vert{}Q^{-1}-\id\Vert{}
  \qquad
  &&\text{($P\in\SO(3)$)}\\
  &\geq
  \left|
   \hinner{\frac{F_{j-1}(p)}{|F_{j-1}(p)|}}{%
     \vect{u}}\right|-|\Theta|
  \geq
  \left|
   \hinner{\frac{F_{j-1}(p)}{|F_{j-1}(p)|}}{%
     \vect{u}}\right|-\frac{4}{\sqrt{N}}\qquad
  &&\text{(Lemma~\ref{lem:rotation})}\\
  &
  \geq
     1-\frac{c_1}{\sqrt{N}}-\frac{4}{\sqrt{N}}
   = 1-\frac{C_2}{\sqrt{N}}.
  &&\text{(Lemma~\ref{lem:axis-vector}, \eqref{eq:C1C2}).}\\
 \end{alignat*}}
 Thus we have the conclusion.
\end{proof}

Now, we apply the ``L\'opez-Ros deformation'' to the holomorphic null
immersion $\widetilde E$.
The following lemma is the straightforward conclusion of
the classical Runge's theorem:
\begin{lemma}\label{lem:runge}
 There exists a holomorphic function $h$ on $\C$
 which does not vanish on $\C$ and 
 satisfies
 \[\left\{
     \begin{array}{rll} 
       |h-1|&\leq \varepsilon_1\qquad
	  &\text{{\rm(}on $\overline{\D}_1\setminus\varpi_j{\rm)}$}\\
       |h-T|&\leq 1\qquad
	  &\text{{\rm(}on $\omega_j${\rm)}}
     \end{array}
   \right.,
 \]
 where
 \begin{equation}\label{eq:runge-constant}
     \varepsilon_1 = 
      \frac{\varepsilon}{\varepsilon+4\sqrt{2}\mu_{j-1} N^{9/4}},\qquad
     \mu_{j-1} =\max_{\overline{\D}_1}|\varphi_{j-1}|,\qquad
     T = 4N^{7/2}+1.
 \end{equation}
\end{lemma}
Using the function $h$ in Lemma~\ref{lem:runge} as a L\'opez-Ros parameter, 
we produce new Weierstrass data as follows:
\begin{equation}\label{eq:runge}
   \hat g := \frac{g}{h},\qquad
   \hat \eta := h\eta,\qquad
   \hat\psi :=
   \frac{1}{2}\bigl(
               1-\hat g^2,\imag(1+\hat g^2), 2\hat g
              \bigr)\hat\eta.
\end{equation}
We denote
\begin{equation}\label{eq:f-j}
    \hat E(z):=\int_0^z \hat\psi(z)\,dz,\qquad
    F_{j}:= B \hat E,
\end{equation}
where $B$ is the matrix as in \eqref{eq:matrix-B}.
By definition \eqref{eq:runge}, $g\eta=\hat g\hat \eta$ holds.
Thus, if we write
\[
   \tilde\psi = (\tilde\psi^{(1)},\tilde\psi^{(2)},\tilde\psi^{(3)})
   \qquad\text{and}\qquad
   \hat\psi = (\hat\psi^{(1)},\hat\psi^{(2)},\hat\psi^{(3)}),
\]
then 
\begin{equation}\label{eq:psi-third}
   \tilde\psi^{(3)}=
   \hat\psi^{(3)}
\end{equation}
holds.

Now, the construction procedure of $F_j$ is accomplished.
Thus, we obtain a sequence $\{F_j\}_{j=0,1,\dots,2N}$
of holomorphic null immersions and a sequence 
$\{\vect{v}_j\}_{j=1,\dots,2N}$ of unit vectors.

In this subsection, we shall prove that $\{F_j\}$ and $\{\vect{v}_j\}$ 
satisfy the conclusions \ref{key:concl:0}--\ref{key:concl:6}
of the Key Lemma~\ref{lem:key}.

\begin{lemma}[\ref{key:concl:6}]\label{lem:c-6}
 For each $j=1,\dots,2N$,
 $\hinner{F_{j}}{\vect{v}_j}=\hinner{F_{j-1}}{\vect{v}_j}$ holds.
\end{lemma} 
\begin{proof}
 By \eqref{eq:psi-third}, we have
 \begin{equation}\label{eq:third}
    \inner{\hat E}{\vect{e}_3}=\int_0^z\hat\psi^{(3)}(w)\,dw=
       \int_0^z\tilde\psi^{(3)}(w)\,dw=
    \inner{\tilde E}{\vect{e}_3}\quad(\vect{e}_3=(0,0,1)).
 \end{equation}
 Since $B\in (e^{-\imag\tau}\delta)\cdot \O(3,\C)$,
 \begin{equation}\label{eq:b-isom}
   \inner{B\vect{x}}{B\vect{y}}=e^{-2\imag\tau}\delta^2\inner{\vect{x}}{\vect{y}}
 \end{equation}
 holds.
 Then
 \begin{alignat*}{2}
  &\hinner{{F}_j}{\vect{v}_j}=
  \inner{{F}_j}{\bar{\vect{v}}_j}\\
  &=
  \inner{B\hat E}{\bar{\vect{v}}_j}
  =\inner{B\hat E}{\vect{b}^{(3)}}
  =\inner{B\hat E}{B\vect{e}_3}\qquad 
  &&\text{(by \eqref{eq:f-j},\eqref{eq:v-j},\eqref{eq:matrix-B})}\\
  &=
  e^{-2\imag\tau}\delta^2\inner{\hat E}{\vect{e}_3}=
  e^{-2\imag\tau}\delta^2\inner{\tilde E}{\vect{e}_3}
    =
  \hinner{F_{j-1}}{\vect{v}_j}
  \quad &&\text{(by \eqref{eq:b-isom},\eqref{eq:third})}.\qed
 \end{alignat*}
\renewcommand{\qed}{\relax}
\end{proof}
The properties \ref{key:concl:4}, \ref{key:concl:5} and
\ref{key:concl:6} in the Key Lemma~\ref{lem:key}
for $l=j$ hold by Lemmas~\ref{lem:c-45} and \ref{lem:c-6}.
The property \ref{key:concl:1} holds trivially because of \eqref{eq:f-j}.
So we shall prove that $F_j$ satisfies \ref{key:concl:2}, \ref{key:concl:3}
of the Key Lemma~\ref{lem:key}.
\begin{lemma}\label{lem:key-c23}
 The holomorphic null immersion $F_j$ as in 
 \eqref{eq:f-j} satisfies
 \begin{alignat}{2}
  \label{eq:fin:1}
        |\varphi_j-\varphi_{j-1}|&\leq\frac{\varepsilon}{2N^2}
        \qquad&&\text{on $\overline{\D}_1\setminus\varpi_j$},\\
  \label{eq:fin:2}
	|\varphi_j|&\geq \frac{C_1}{N^{3/4}}
        \qquad&&\text{on $\overline{\varpi}_j$},\\
  \label{eq:fin:3}
        |\varphi_j|&\geq C_1N^{9/4}
        \qquad&&\text{on $\omega_j$}.
 \end{alignat}
where $C_1$ is given in \eqref{eq:C1C2}.
\end{lemma}
\begin{proof}
 By the definitions \eqref{eq:cal-f}, 
 \eqref{eq:tildeE} and \eqref{eq:f-j},
 and noticing that $Q^{-1}$ and $A$ commute (cf.\ \eqref{eq:commute}),
 we have
 \begin{align*}
    \varphi_{j-1}&=e^{-\imag\tau}P^{-1}AQ^{-1}\tilde\psi=
               e^{-\imag\tau}P^{-1}Q^{-1}A\tilde\psi
                  ,\\
    \varphi_{j}&=e^{-\imag\tau}P^{-1}AQ^{-1}\hat\psi=
                   e^{-\imag\tau}P^{-1}Q^{-1}A\hat\psi.
 \end{align*}
 Then 
 \begin{equation}\label{eq:rel-norm}
   |\varphi_{j-1}|=|A\tilde\psi|,\quad |\varphi_{j}|=|A\hat\psi|,
   \quad 
 |\tilde\psi|=  |A^{-1}P\varphi_{j-1}|,\quad |\hat\psi|=|A^{-1}P\varphi_{j}|,
 \end{equation}
 hold because 
 $P$, $Q\in\SO(3)$.
 By \eqref{eq:norm-ineq} and \eqref{eq:a-norm}, 
 \begin{equation}\label{eq:rel-norm-2}
  \begin{aligned}
  |\varphi_{j}-\varphi_{j-1}|&=
  |A(\hat\psi-\tilde\psi)|\leq \Vert{}A\Vert{}\,|\hat\psi-\tilde\psi|
  \leq \sqrt{2}|\hat\psi-\tilde\psi|\\
  |\varphi_{j}-\varphi_{j-1}|&=
  |A(\hat\psi-\tilde\psi)|\geq \frac{1}{\Vert{}A^{-1}\Vert{}}|\hat\psi-\tilde\psi|
  \geq \frac{1}{\sqrt{2}{{N^{1/4}}}}|\hat\psi-\tilde\psi|
  \end{aligned}
 \end{equation}
 hold.
 Here, by \eqref{eq:weierstrass-f}, \eqref{eq:runge} and 
\eqref{eq:psi-third},
 we have
 \begin{alignat*}{2}
  |&\hat\psi-\tilde\psi|
   =\left|
     \frac{1}{2}\biggl((1-\hat g^2)\hat\eta - (1-g^2)\eta,
                      \imag(1+ \hat g^2)\hat \eta-\imag(1+g^2)\eta\biggr)
    \right|\\
    &=\frac{1}{2}
    \left|
            \left(\left(1-\frac{g^2}{h^2}\right)h\eta - (1-g^2)\eta,
             \imag\left(1+ \frac{g^2}{h^2}\right)h\eta-\imag(1+g^2)\eta\right)
    \right|
    \\
   &=\frac{1}{2}
      \left|
        (h-1)
        \left(
          \left(1+\frac{g^2}{h}\right),
          \imag\left(1-\frac{g^2}{h}\right)
        \right)\eta
      \right|
  \\
   &= \frac{1}{2}|h-1|\,|\eta|\,
          \left(
          \left|1+\frac{g^2}{h}\right|^2+
          \left|1-\frac{g^2}{h}\right|^2
          \right)^{1/2}\!\!
   \leq 
          \frac{1}{2}|h-1|\,|\eta|\,
          \left(
          \left|1+\frac{g^2}{h}\right|+
          \left|1-\frac{g^2}{h}\right|
          \right)\\
   & \leq  |h-1|\,|\eta|\,\left(1+\frac{|g|^2}{|h|}\right)
    \leq
     |h-1|\,|\eta|
     \left(
        1+ \frac{|g|^2}{1-|h-1|}
     \right)\\
     &\leq 
      |h-1|
      \frac{(1+|g|^2)|\eta|}{1-|h-1|}=\sqrt{2}|\tilde\psi|\frac{|h-1|}{1-|h-1|}.
 \end{alignat*}
 Since $h$ is taken as in Lemma~\ref{lem:runge} and $P\in \rm{SO(3)}$, we have
 \begin{alignat*}{2}
  |\hat\psi-\tilde\psi|
   &\leq\sqrt{2}|\tilde\psi|\,\frac{\varepsilon_1}{1-\varepsilon_1}
   =\sqrt{2}|\tilde\psi|\,\frac{\varepsilon}{4\sqrt{2}
           \mu_{j-1} N^{9/4}}\quad
   &&\text{(Lemma~\ref{lem:runge}, \eqref{eq:runge-constant})}\\
   & =|A^{-1}P \varphi_{j-1}|\,\frac{\varepsilon}{4\mu_{j-1} N^{9/4}}
    \leq\Vert{}A^{-1}\Vert{}\,|P \varphi_{j-1}|\,
     \frac{\varepsilon}{4\mu_{j-1} N^{9/4}}
  \quad &&\text{(by \eqref{eq:rel-norm},\eqref{eq:norm-ineq})}\\
   & \leq {\sqrt{2}{{N^{1/4}}}}\,|\varphi_{j-1}|
      \,\frac{\varepsilon}{4\mu_{j-1} N^{9/4}}
  \qquad &&\text{(by \eqref{eq:a-norm})}\\
   &\leq \sqrt{2}{{N^{1/4}}}\mu_{j-1}\frac{\varepsilon}{4\mu_{j-1} N^{9/4}}
    = \frac{\varepsilon}{2\sqrt{2}N^2}
  \qquad &&\text{(by \eqref{eq:runge-constant})}
 \end{alignat*}
 holds on $\overline{\D}_1\setminus\varpi_j$.
 Thus, by \eqref{eq:rel-norm-2}
 $|\varphi_j-\varphi_{j-1}|
   \leq \sqrt{2}|\hat\psi-\tilde\psi|
   \leq {\varepsilon}/{(2N^2)}$,
 which is \eqref{eq:fin:1}.

 Next, on $\overline{\varpi}_j$, 
 it holds that
\allowdisplaybreaks{%
 \begin{alignat*}{2}
   |\varphi_j|&=|A\hat\psi|
       \geq \frac{1}{\Vert{}A^{-1}\Vert{}}|\hat\psi|
       \geq \frac{1}{\sqrt{2}{{N^{1/4}}}}|\hat\psi|
     \qquad&&
  \text{(by \eqref{eq:rel-norm},\eqref{eq:norm-ineq},\eqref{eq:a-norm})}\\
     & = \frac{1}{\sqrt{2}{{N^{1/4}}}}
                   \frac{1}{\sqrt{2}}(1+|\hat g|^2)|\,|\hat\eta|
      = \frac{1}{2{{N^{1/4}}}}
                   (1+|\hat g|^2)\,|\hat\eta|
        \qquad&&\text{(by \eqref{eq:runge}, \eqref{eq:induced-weier})}\\
     & \geq \frac{1}{{{N^{1/4}}}}|\hat g \hat\eta|
                    = \frac{1}{{{N^{1/4}}}}|g\eta|
    \qquad&&\text{(by \eqref{eq:psi-third})}
  \\
            &= \frac{\sqrt{2}}{{{N^{1/4}}}}\frac{1}{\sqrt{2}}(1+|g|^2)|\eta|
                  \frac{|g|}{1+|g|^2}
            = \frac{\sqrt{2}|\tilde\psi|}{{{N^{1/4}}}}\frac{|g|}{1+|g|^2}
  \qquad&&\text{(by \eqref{eq:weierstrass-f})}\\
            & =
  \frac{\sqrt{2}|A^{-1} P \varphi_{j-1}|}{{{N^{1/4}}}}\frac{|g|}{1+|g|^2}
  \geq \frac{\sqrt{2}|\varphi_{j-1}|}{{{N^{1/4}}}\Vert{}A\Vert{}}
       \frac{|g|}{1+|g|^2}
  && \! \!\!   \! \!\!   \! \!\!   \text{(by \eqref{eq:rel-norm}, \eqref{eq:norm-ineq}, $P\in \rm{SO(3)}$)}\\
  &\geq \frac{|\varphi_{j-1}|}{{{N^{1/4}}}}
                  \frac{|g|}{1+|g|^2}
  \geq
               \frac{|\varphi_{j-1}|}{{{N^{1/4}}}}\frac{\frac{1}{2\sqrt{N}}}{%
                        1+\frac{1}{4N}}
  &&\quad \text{(by \eqref{eq:a-norm}, Lemma~\ref{lem:gauss-bound})}\\
  &\geq
  \frac{\nu}{2{{N^{1/4}}}}\frac{\frac{1}{2\sqrt{N}}}{%
  1+\frac{1}{4N}}
  =\frac{\nu}{N^{3/4}}\frac{1}{4+1/N}\geq \frac{\nu}{5N^{3/4}}=
          \frac{C_1}{N^{3/4}}
  \quad&& \quad  \quad\text{(Lemma~\ref{lem:bound-j}, \eqref{eq:C1C2})}.
 \end{alignat*}}%
 Thus, we have \eqref{eq:fin:2}.

 Finally, on  $\omega_j$, we have
 \[
     |\varphi_j|\geq \frac{1}{2{{N^{1/4}}}}(1+|\hat g|^2)|\hat\eta|
 \]
 by the same way as in the previous argument.
 Then
 \allowdisplaybreaks{
 \begin{alignat*}{2}
   |\varphi_j|&\geq \frac{1}{2{{N^{1/4}}}} (1+|\hat g|^2)|\hat\eta|\geq
               \frac{1}{2{{N^{1/4}}}}|\hat\eta|\\
              &=\frac{1}{2{{N^{1/4}}}}|h|\,|\eta|
               =\frac{\sqrt{2}}{2{{N^{1/4}}}}
               \frac{1}{\sqrt{2}}(1+|g|^2)|\eta|\frac{|h|}{1+|g|^2}
          && \quad \qquad     \qquad \text{(by \eqref{eq:runge})}\\ 
              & =\frac{|h|}{\sqrt{2}{{N^{1/4}}}}
               |\tilde\psi|\frac{1}{1+|g|^2} 
               = \frac{|h|}{\sqrt{2}{{N^{1/4}}}}
                |A^{-1}P \varphi_{j-1}|\frac{1}{1+|g|^2}
      &&   \qquad \text{(by \eqref{eq:weierstrass-f}, \eqref{eq:rel-norm})}\\ 
              &\geq
               \frac{|h|}{\Vert{}A\Vert{}\sqrt{2}{{N^{1/4}}}}|P \varphi_{j-1}|\frac{1}{1+|g|^2}\geq 
               \frac{|h|}{2{{N^{1/4}}}}|P\varphi_{j-1}|\frac{1}{1+|g|^2}
             &&   \,\, \, \quad   \quad\text{(by \eqref{eq:norm-ineq}, \eqref{eq:a-norm})}\\ 
              &\geq
               \frac{|h|}{2{{N^{1/4}}}}|\varphi_{j-1}|\frac{1}{1+4N} 
               \geq
               \frac{|h|}{2{{N^{1/4}}}}\frac{\nu}{2}\frac{1}{1+4N}
            &&\!\!\!\! \! \!\!\!\! \!\!\!\!\! \!  \!\!\!\! \text{($P\in \rm{SO(3)}$,   Lemmas~\ref{lem:gauss-bound}, \ref{lem:bound-j})}\\ 
             &\geq\frac{|h|\nu}{4{{N^{1/4}}}}\frac{1}{5N}=\frac{\nu}{20N^{5/4}}|h|
                     =\frac{\nu}{20N^{5/4}}(T-|h-T|)
             \\
             &\geq\frac{\nu}{20N^{5/4}}4N^{7/2}=\frac{\nu}{5}N^{9/4}
             =C_1N^{9/4}
 && \quad \text{(Lemma~\ref{lem:runge}, \eqref{eq:C1C2})}.\\
 \end{alignat*}}%
 Hence we have \eqref{eq:fin:3}.
\end{proof}
Thus we have $\{F_j\}$ and $\{\vect{v}_j\}$ satisfying
properties \ref{key:concl:0}--\ref{key:concl:6}
in Lemma~\ref{lem:key}.

\section{%
 A proof of Proposition~\ref{prop:main}}
\label{sec:proof-main}
In this section, we prove Proposition~\ref{prop:main}.
We take the sequences $\{F_j\}$ and $\{\vect{v}_j\}$ as in the Key
Lemma~\ref{lem:key}, and set 
\begin{equation}\label{eq:Y}
   Y := F_{2N}.
\end{equation}
Recall that $X =  F_0$  by \ref{key:concl:0}. Then we shall prove 
\ref{key:concl:d1}--\ref{key:concl:d3} in Proposition~\ref{prop:main}.
\begin{lemma}\label{lem:d1}
 It holds that 
 \[
    |\varphi_{Y}-\varphi_X|\leq 
    \frac{\varepsilon}{N},
    \quad\text{and}\quad
    |Y-X|\leq \frac{2\varepsilon}{N}
    \qquad
    \text{on 
    $\overline{\D}_1\setminus\bigl(\varpi_1\cup\dots\cup\varpi_{2N}\bigr)$}.
 \]
\end{lemma}
\begin{proof}
 By \ref{key:concl:2} of the Key Lemma~\ref{lem:key}, 
 \[
   |\varphi_Y-\varphi_X|=
   |\varphi_{2N}-\varphi_0|
   \leq |\varphi_{2N}-\varphi_{2N-1}|+\dots+|\varphi_1-\varphi_0|
   \leq 2N\cdot \frac{\varepsilon}{2N^2}=\frac{\varepsilon}{N}
 \]
 holds on $\overline{\D}_1\setminus\bigl(\varpi_1\cup\dots\cup\varpi_{2N}\bigr)$.
 On the other hand, by Lemma~\ref{lem:fj-est},
 \[
   |Y-X|=
   |F_{2N}-F_0|
   \leq |F_{2N}-F_{2N-1}|+\dots+|F_1-F_0|
   \leq 2N\cdot \frac{\varepsilon}{N^2}=\frac{2\varepsilon}{N}
 \]
 holds on $\overline{\D}_1\setminus\bigl(\varpi_1\cup\dots\cup\varpi_{2N}\bigr)$.
\end{proof}
\begin{corollary}[the conclusion \ref{key:concl:d1}]\label{cor:d1}
 It holds that
 \[
     |\varphi_Y-\varphi_X|<\varepsilon
      \qquad\text{and}\qquad
     |Y-X|<\varepsilon
     \qquad \text{on  \quad $\D_{1-\varepsilon}$}.
 \]
\end{corollary}
\begin{proof}
 Note that we take the labyrinth as in Appendix~\ref{app:labyrinth}.
 Here, by \eqref{eq:N3}, 
 \[
     \frac{2}{N}+\frac{1}{8N^3}=
     \frac{1}{N}\left(2+\frac{1}{8N^2}\right)
     < \frac{3}{N}\leq\varepsilon
 \]
 holds.
 Then by \ref{item:omega:2} of Lemma~\ref{lem:omega} in
 Appendix~\ref{app:labyrinth}, we have that
 \begin{equation}\label{eq:d-1-e}
    \D_{1-\varepsilon}\subset
    \overline{\D}_1\setminus
           \bigl(
            \varpi_1\cup\dots\cup\varpi_{2N}
           \bigr).
 \end{equation}
 Thus, by Lemma~\ref{lem:d1} and \eqref{eq:N1}, 
 it holds on $\D_{1-\varepsilon}$ that 
 \[
    |\varphi_Y-\varphi_X|
    =|\varphi_{2N}-\varphi_0|\leq 
       \frac{\varepsilon}{N}\leq\varepsilon,\qquad
    |Y-X|
    =|F_{2N}-F_0|\leq 
       \frac{2\varepsilon}{N}\leq\varepsilon
   \qquad\qed
 \]
 \renewcommand{\qed}{\relax}
\end{proof}


\begin{lemma}\label{lem:phi-y-est}
 The function $\varphi_Y=\varphi_{2N}$ satisfies
 \[
     |\varphi_Y|\geq 
     \begin{cases}
        \dfrac{C_1}{2}N^{9/4}\qquad &
          \text{on $\omega_1\cup\dots\cup\omega_{2N}$}\\[10pt]
        \dfrac{C_1}{2N^{3/4}}\qquad &
          \text{on $\overline{\D}_1$}.
      \end{cases}
 \]
\end{lemma}
\begin{proof}
 On $\omega_j$,
 \begin{alignat*}{2}
  |\varphi_{Y}|&=|\varphi_{2N}|
      \geq |\varphi_j|-
           |\varphi_{2N}-\varphi_{2N-1}|-\dots-
           |\varphi_{j+1}-\varphi_j|&&\\
     &\geq C_1N^{9/4} - \frac{(2N-j+1)\varepsilon}{2N^2}&&
          (\text{by \ref{key:concl:3}, \ref{key:concl:2}})\\
     &\geq C_1N^{9/4}-\frac{\varepsilon}{N}
       =N^{9/4} \left(C_1-\frac{\varepsilon}{N^{3/4}}\right)&&\\
     &\geq 
        N^{9/4} \left(C_1-\frac{\varepsilon}{N^{1/4}}\right)\geq 
         \frac{C_1}{2}N^{9/4}
      \qquad &&\text{(by \eqref{eq:N3})}.
 \end{alignat*}
 On the other hand, on $\varpi_j$, we have
 \begin{align*}
  |\varphi_{Y}|&=|\varphi_{2N}|\geq 
  |\varphi_j|-|\varphi_{2N}-\varphi_{2N-1}|-\dots-|\varphi_{j+1}-\varphi_j|
    &&\\
   &\geq \frac{C_1}{N^{3/4}} - \frac{(2N-j+1)\varepsilon}{2N^2}&&
          (\text{by \ref{key:concl:3}, \ref{key:concl:2}})\\
   &\geq \frac{C_1}{N^{3/4}}-\frac{\varepsilon}{N}
   =\frac{1}{N^{3/4}} \left(C_1-\frac{\varepsilon}{N^{1/4}}\right)
    \geq \frac{C_1}{2N^{3/4}} \qquad &&\text{(by \eqref{eq:N3})}.
 \end{align*}
 Finally, on  
 $\overline{\D}_1\setminus(\varpi_1\cup\dots\cup\varpi_{2N})$,%
 \begin{alignat*}{2}
   |\varphi_Y|&=|\varphi_{2N}|\geq
   |\varphi_0|-|\varphi_{2N}-\varphi_{2N-1}|-\dots-|\varphi_{1}-\varphi_0|&&
   \\
   &\geq |\varphi_0|-2N\cdot\frac{\varepsilon}{2N^2}
   \geq \nu - \frac{\varepsilon}{N}
     \qquad&&\text{(\ref{key:concl:2}, \eqref{eq:constants})}\\
   &\geq 
         \nu-\frac{\varepsilon}{N^{3/4}}
   \geq \frac{C_1}{2N^{3/4}}
     \qquad&&\text{(by \eqref{eq:N3})}.
 \end{alignat*}
 Hence we have the conclusion.
\end{proof}
\begin{corollary}[the conclusion \ref{key:concl:d2}]
\label{cor:intrinsic}
 The disc
 $(\D_1,ds^2_Y)$ contains 
 a geodesic disc $\DD$ centered at $0$ with radius $\rho+s$.
\end{corollary}
\begin{proof}
 The induced metric $ds^2_Y$ is expressed as 
 \[
   ds^2_{Y}=|\varphi_{Y}|^2\,|dz|^2.
 \]
 Consider a Riemannian metric
 \[
    ds^2 := \left(\frac{2N^{3/4}}{C_1}\right)^2 ds^2_Y
          = \lambda^2 \,|dz|^2,\qquad
      \left(\lambda:= \frac{2N^{3/4}}{C_1}|\varphi_{Y}|\right).
 \] 
 Then by Lemma~\ref{lem:phi-y-est},
 $ds^2$ satisfies the assumptions of Lemma~\ref{lem:length-estimate} 
 in Appendix~\ref{app:labyrinth}.
 Thus, we have
 \[
    \dist_{ds^2}(0,\partial\overline{\D}_1)\geq N,
 \]
 where $\dist_{ds^2}$ denotes the distance function 
 with respect to $ds^2$.
 Then by \eqref{eq:N3}, we have
 \[
   \dist_{ds^2_{2N}}(0,\partial\overline{\D}_1)\geq \frac{C_1}{2N^{3/4}}N
    =\frac{C_1N^{1/4}}{2}
   \geq \rho+s.
 \]
 Hence we have the conclusion.
\end{proof}

\begin{figure}
   \begin{center}
    \includegraphics[width=0.6\textwidth]{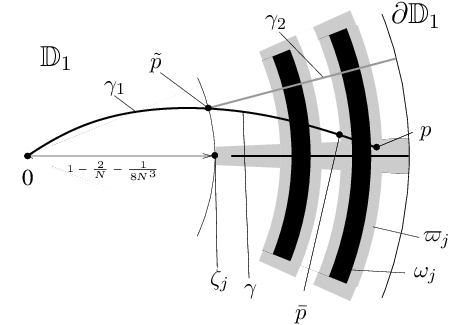}
   \end{center}
 \caption{The curve $\gamma$ and points $\bar p$, $\tilde p$}
 \label{fig:curve-gamma}
\end{figure}

By Corollary~\ref{cor:intrinsic}, 
one can take a geodesic disc $\DD$ of $(\D_1,ds^2_{Y})$
centered at the origin with radius $\rho+s$.
Fix $p\in\partial\DD$, and prove \ref{key:concl:d3} of the Key
Lemma~\ref{lem:key}.
First, we assume $p\in\varpi_j$ for some $j\in \{1,\cdots,2N\}$
(otherwise, the proof of \ref{key:concl:d3} is rather easy).
Since $p\in\partial\DD$,
there exists a $ds^2_Y$-geodesic $\gamma$ joining $0$ and $p$
with length $\rho+s$.
Since $ds^2_{Y}$ is a Riemannian metric of non-positive
Gaussian curvature, 
\begin{equation}\label{eq:shortest}
\text{
\begin{minipage}{0.8\textwidth}
 an arbitrary subarc of $\gamma$ 
  is the shortest geodesic.
\end{minipage}}
\end{equation}
Hence the image of $\gamma$ is contained in $\DD$.
\begin{lemma}\label{lem:gamma-length}
 The Euclidean length of  $\gamma$ satisfies
 \[
     \Length_{\C}(\gamma)\leq \frac{2(\rho+s)}{C_1} N^{3/4}.
 \]
\end{lemma}
\begin{proof}
  Since he $ds^2_{Y}$-arclength of $\gamma$ is  $\rho+s$,
  Lemma \ref{lem:phi-y-est} implies that
 \[
     \rho+s = \int_{\gamma}|\varphi_{Y}|\,|dz|
            \geq \int_{\gamma}\frac{C_1}{2N^{3/4}}\,|dz|
            =\frac{C_1}{2N^{3/4}}\Length_{\C}(\gamma).
 \]
  Hence we have the conclusion.
\end{proof}
Now, take points  $\bar p$, $\tilde p\in\DD$ on the arc $\gamma$
such that
\begin{itemize}
 \item $\bar p\in\partial\varpi_j$ and 
       the subarc of $\gamma$ joining $\bar p$ and $p$ 
       is contained in $\overline{\varpi}_j$,
       namely, $\bar p$ is the final point where $\gamma$
       meets $\partial \varpi_j$, 
 \item and the subarc of $\gamma$ joining $0$ and 
       $\tilde p\in\partial\D_{1-\frac{2}{N}-\frac{1}{8N^3}}$
       contained in $\overline{\D}_{1-\frac{2}{N}-\frac{1}{8N^3}}$;
       namely, $\tilde p$ is the first point where 
       $\gamma$ meets $\partial\D_{1-\frac{2}{N}-\frac{1}{8N^3}}$.
\end{itemize}
See Figure~\ref{fig:curve-gamma}.
\begin{lemma}\label{lem:length}
 It holds that 
 \begin{align}
  \label{eq:length-1}
    &|F_{l}(\bar p)|\leq r + \frac{2\varepsilon}{N}
    \qquad (l=0,\dots,2N),\\
      \label{eq:length-1bis}
    &|F_{j-1}( p)|\leq r + \frac{2\varepsilon}{N}
  ,\\
    \label{eq:length-2}
    &|F_{2N}(p)-F_{j}(p)|\leq \frac{2\varepsilon}{N},\qquad
    |F_{2N}(\bar p)-F_{j}(\bar p)|\leq \frac{2\varepsilon}{N},\\
  \label{eq:length-3}
    &|F_{2N}(p)-F_{2N}(\bar p)|\leq
    s + \frac{c_2}{N^{1/4}},
 \end{align}
 where 
 $C_1$ and $c_2$ are defined by \eqref{eq:C1C2} and \eqref{eq:c123},
 respectively.
\end{lemma}
\begin{proof}
 Since $\bar p\not\in\varpi_1\cup\dots\cup\varpi_{2N}$,
 Lemma~\ref{lem:fj-est} and the assumption \ref{key:ass:3}
 of the Proposition~\ref{prop:main}imply
 \begin{align*}
    |F_l(\bar p) |&\leq
    |F_0(\bar p)|+|F_1(\bar p)-F_0(\bar p)|+\dots+
                  |F_l(\bar p)-F_{l-1}(\bar p)|\\
    &\leq 
      r + \frac{l\varepsilon}{N^2}
      \leq  r + \frac{2\varepsilon}{N}.
 \end{align*}
 Hence we have \eqref{eq:length-1}. A similar reasoning proves  \eqref{eq:length-1bis} .

 Since  $p\not\in\varpi_{j+1}\cup\dots\cup\varpi_{2N}$,
 Lemma~\ref{lem:fj-est} implies
 \begin{align*}
     |F_{2N}(p)-F_{j}(p)|&\leq |F_{2N}(p)-F_{2N-1}(p)|+\dots+
        |F_{j+1}(p)-F_{j}(p)|\\
     &\leq \frac{(2N-j)\varepsilon}{N^2}<\frac{2\varepsilon}{N}.
 \end{align*}
 Then the first inequality of \eqref{eq:length-2} holds.
 Similarly, we have the second inequality of \eqref{eq:length-2}.

 Let $\gamma_1$ be the subarc of the geodesic $\gamma$ joining $0$
 and $\tilde p$, and let $\gamma_2$ be the line segment joining $\tilde p$
 and $\partial\D_1$ which is contained in the line $0\tilde p$,
 see Figure~\ref{fig:curve-gamma}.
 Since $\gamma_1\cup\gamma_2$ is a path joining $0$ and $\partial\D_1$,
 the assumption \ref{key:ass:2} and \ref{key:concl:0} implies that
 \begin{equation}\label{eq:init-dist}
     \Length_{ds_0^2}(\gamma_1\cup\gamma_2)
     =\int_{\gamma_1\cup\gamma_2}|\varphi_0(z)|\,|dz|
     \geq \dist_{ds_0^2}(0,\partial\D_1)
     \geq \rho,
 \end{equation}
where  $\Length_{ds_0^2}(\gamma_1\cup\gamma_2)$ is the
length of the curve $\gamma_1\cup\gamma_2$ 
with respect to the metric $ds^2_0$.
 On the other hand, by \eqref{eq:constants}, we have
 \begin{align}\label{eq:gamma-2}
     \Length_{ds_0^2}(\gamma_2)&=
     \int_{\gamma_2}|\varphi_0(z)|\,|dz|
     \leq \mu\cdot\Length_{\C}(\gamma_2)\\
     &= \mu\left(\frac{2}{N}+\frac{1}{8N^3}\right)
     =\frac{\mu}{N}\left(2+\frac{1}{8N^2}\right)
     \leq \frac{3\mu}{N}.\nonumber
 \end{align}
 Hence we have 
 \begin{equation}\label{eq:phi-0-gamma-1}
  \int_{\gamma_1}|\varphi_0(z)|\,|dz|
   =
  \int_{\gamma_1\cup\gamma_2}|\varphi_0(z)|\,|dz|-
       \int_{\gamma_2}|\varphi_0(z)|\,|dz|
   \geq \rho-\frac{3\mu}{N}.
 \end{equation}
 Since $\gamma_1$ is contained in the subarc
 of $\gamma$ joining $0$ and $\bar p$, 
 we have
 \begin{alignat*}{2}
  \dist_{ds^2_{2N}}(0,\bar p) 
  &\geq 
  \dist_{ds^2_{2N}}(0,\tilde p) =
  \int_{\gamma_1}|\varphi_{2N}(z)|\,|dz|
  &&\text{(by \eqref{eq:shortest})}\\
  &=
  \int_{\gamma_1}(|\varphi_{2N}(z)|-|\varphi_0(z)|)\,|dz|+
  \int_{\gamma_1}|\varphi_0(z)|\,|dz|\\
  &\geq 
  -\int_{\gamma_1}|\varphi_{2N}(z)-\varphi_0(z)|\,|dz|+
  \int_{\gamma_1}|\varphi_0(z)|\,|dz|\\
  &\geq 
  -\int_{\gamma_1}|\varphi_{2N}(z)-\varphi_0(z)|\,|dz|+
      \rho-\frac{3\mu}{N}
  \qquad&&\text{(by \eqref{eq:phi-0-gamma-1})}
  \\
  &\geq
  -\int_{\gamma_1}\frac{\varepsilon}{N}\,|dz|+\rho-\frac{3\mu}{N}
  \qquad &&\text{(Lemma~\ref{lem:d1})}\\
  &\geq 
   -\Length_{\C}(\gamma)\frac{\varepsilon}{N}+\rho-\frac{3\mu}{N}
  \qquad&&\text{($\gamma_1\subset\gamma$)}\\
  &\geq 
    -\frac{2\varepsilon (\rho+s)}{C_1 N^{1/4}}+\rho -\frac{3\mu}{N}
  \geq
    -\frac{2\varepsilon (\rho+s)}{C_1 N^{1/4}}+\rho -\frac{3\mu}{N^{1/4}}
  \qquad &&\text{(Lemma~\ref{lem:gamma-length})}\\
  &= \rho -\frac{1}{N^{1/4}}\left(3\mu+\frac{2\varepsilon(\rho+s)}{C_1}\right)
   = \rho-\frac{c_2}{N^{1/4}}.\qquad
  &&\text{(by \eqref{eq:c123})}
 \end{alignat*}
 Here, since $\bar p$ lies on the geodesic $\gamma$ joining $0$ and $p$,
 \eqref{eq:shortest} implies
 \begin{align*}
  |F_{2N}(p)-F_{2N}(\bar p)|&\leq 
  \dist_{ds^2_{2N}}(p,\bar p) =
  \dist_{ds^2_{2N}}(0,p)-\dist_{ds^2_{2N}}(0,\bar p)\\
   &= \rho +s - \dist_{ds^2_{2N}}(0,\bar p)\\
   &\leq \rho+s -\left(\rho-\frac{c_2}{N^{1/4}}\right)
   = s+\frac{c_2}{N^{1/4}}.
 \end{align*}
 Thus \eqref{eq:length-3} is obtained.
\end{proof}

\subsubsection*{Case 1: $p\in\varpi_j$ and $|F_{j-1}(\zeta_j)|>1/\sqrt{N}$}
\begin{lemma}\label{lem:d3-case-1}
 When $p\in\varpi_j$ and $|F_{j-1}(\zeta_j)|>1/\sqrt{N}$,
 \[
     |F_j(p)|\leq \sqrt{r^2+s^2}+\frac{c_3}{{{N^{1/4}}}}
 \]
 holds.
\end{lemma}
\begin{proof}
 Let $\vect{v}_j\in\C^3$ be the unit vector
 as in \ref{key:concl:4}--\ref{key:concl:6}, and denote
 \[
    (\vect{v}_j)^{\perp}:=
    \text{%
         (the orthogonal complement of 
          $\vect{v}_j$ with respect to $\hinner{~}{~}$)}.
 \]
 Then $(\vect{v}_j)^{\perp}$ is a (complex) $2$-dimensional subspace
 of $\C^3$.
 Denote by $\Pi_j$ the orthogonal projection 
 \begin{equation}\label{eq:proj}
    \Pi_j\colon{}\C^3\ni \vect{x}\longmapsto 
    \vect{x}-\hinner{\vect{x}}{\vect{v}_j}\vect{v}_j\in (\vect{v}_j)^{\perp}
 \end{equation}
 with respect to the Hermitian inner product $\hinner{~}{~}$.
 Then for any vector $\vect{x}\in\C^3$, 
 \begin{equation}\label{eq:pythagorean}
  |\vect{x}|^2 = |\hinner{\vect{x}}{\vect{v}_j}|^2 +
                 |\Pi_j\vect{x}|^2
 \end{equation}
 holds.
 Thus, we have
\allowdisplaybreaks{%
 \begin{alignat*}{2}
  \left|
   \Pi_j F_j(p)
  \right|&\leq 
  \left|
   \Pi_j F_j(p)-\Pi_j F_j(\bar p)\right|+
   \left|\Pi_j F_j(\bar p)\right|\\
  &\leq
  \left|
   \Pi_j \bigl(F_j(p)-F_j(\bar p)\bigr)  \right|+
   \left|\Pi_j F_j(\bar p)\right|\\%
  &\leq 
  \left|
   F_j(p)-F_j(\bar p)
  \right|+\left|\Pi_j F_j(\bar p)\right|
  \qquad 
  &&\text{(by \eqref{eq:pythagorean})}\\
  &\leq 
  \left|
   F_j(p)-F_j(\bar p)
  \right|+\left|\Pi_j \bigl(F_{j}(\bar p)-F_{j-1}(\bar p)\bigr)\right|
    +|\Pi_j F_{j-1}(\bar p)|\\
  &\leq 
  \left|
   F_j(p)-F_j(\bar p)
  \right|+\left|F_{j}(\bar p)-F_{j-1}(\bar p)\right|
    +|\Pi_j F_{j-1}(\bar p)|
  \qquad 
  &&\text{(by \eqref{eq:pythagorean})}\\
  &\leq 
  \left|
   F_{2N}(p)-F_{2N}(\bar p)
  \right|+
  \left|
   F_{2N}(p)-F_{j}(p)
  \right|+
  \left|
   F_{2N}(\bar p)-F_{j}(\bar p)
  \right|\\
  &\hphantom{\leq}+
  \left|F_{j}(\bar p)-F_{j-1}(\bar p)\right|
    +|\Pi_j F_{j-1}(\bar p)|\\
  &\leq 
    \left(s+\frac{c_2}{{{N^{1/4}}}}\right)+
    \frac{2\varepsilon}{N}+
    \frac{2\varepsilon}{N}\\
  &\hphantom{\leq++++}+
    \left|F_{j}(\bar p)-F_{j-1}(\bar p)\right|     
    +|\Pi_j F_{j-1}(\bar p)|
   \quad &&\text{(Lemma~\ref{lem:length})}\\
  &\leq 
    \left(s+\frac{c_2}{{{N^{1/4}}}}\right)+
    \frac{4\varepsilon}{N}+
    \frac{\varepsilon}{N^2}
    +|\Pi_j F_{j-1}(\bar p)|
   \quad &&\text{(Lemma~\ref{lem:fj-est})}\\
  &\leq s + \frac{1}{{{N^{1/4}}}}
    \left({c_2}+\frac{4\varepsilon}{N^{3/4}}+
          \frac{\varepsilon}{N^{7/4}}\right)
          +|\Pi_j F_{j-1}(\bar p)|\\
  &\leq s + \frac{c_2+5\varepsilon}{{{N^{1/4}}}}
          +|\Pi_j F_{j-1}(\bar p)|.
\end{alignat*}}
Hence we have
\begin{equation}\label{eq:proj-1}
  \left|
   \Pi_j F_j(p)
  \right|\leq 
    s + \frac{c_2+5\varepsilon}{{{N^{1/4}}}}+|\Pi_j F_{j-1}(\bar p)|.
\end{equation}
Here, we assume $F_{j-1}(\bar p)\neq 0$.
Since $\bar p\in \overline{\varpi}_j$, we have
\begin{alignat*}{2}
 |\Pi_j F_{j-1}(\bar p)|&=
  \sqrt{
  |F_{j-1}(\bar p)|^2 -|\inner{F_{j-1}(\bar p)}{\vect{v}_j}|^2}
  \qquad &&\text{(by \eqref{eq:pythagorean})}\\
 &=|F_{j-1}(\bar p)|
    \sqrt{1-\left|
       \hinner{\frac{F_{j-1}(\bar p)}{|F_{j-1}(\bar p)|}}{\vect{v}_j}
          \right|^2}
 \\
 & \leq 
    |F_{j-1}(\bar p)|
    \sqrt{1-\left(1-\frac{C_2}{\sqrt{N}}\right)^2}
 \qquad && \text{(by \ref{key:concl:5})}\\
  &=
    |F_{j-1}(\bar p)|
    \sqrt{\frac{2C_2}{\sqrt{N}}-\frac{C_2^2}{N}}
   \leq 
    |F_{j-1}(\bar p)|
    \sqrt{\frac{2C_2}{\sqrt{N}}}\\
 &=
    |F_{j-1}(\bar p)|
    \frac{\sqrt{2C_2}}{{{N^{1/4}}}}
    \leq
    \left(r+\frac{2\varepsilon}{N}\right)
      \cdot\frac{\sqrt{2C_2}}{{{N^{1/4}}}}
    &&\text{(by \eqref{eq:length-1} in Lemma~\ref{lem:length})}
    \\
    &\leq 
      \frac{(r+2\varepsilon)\sqrt{2C_2}}{{{N^{1/4}}}}.
 \end{alignat*}
 Then by \eqref{eq:proj-1}, we have
 \begin{equation}\label{eq:fp-proj}
 |\Pi_j F_j(p)|\leq
     s + \frac{\alpha}{{{N^{1/4}}}}
     \qquad \biggl(\alpha:=c_2+5\varepsilon+(r+2\varepsilon)\sqrt{2C_2}\biggr)
 \end{equation}
 when $F_{j-1}(\bar p)\neq 0$.
 Otherwise, namely when $F_{j-1}(\bar p)=0$, 
 \eqref{eq:fp-proj} holds trivially.

 Thus, 
 \begin{alignat*}{2}
  |F_j(p)|&=
  \sqrt{|\hinner{F_j(p)}{\vect{v}_j}|^2+|\Pi_j F_j(p)|^2}
  \qquad &&\text{(by \eqref{eq:pythagorean})}\\
  &=
  \sqrt{|\hinner{F_{j-1}(p)}{\vect{v}_j}|^2+|\Pi_j F_j(p)|^2}
  \qquad &&\text{(by \ref{key:concl:6})}\\
  &\leq 
  \sqrt{|F_{j-1}(p)|^2+|\Pi_j F_j(p)|^2}\\
  &\leq
  \sqrt{\left(r+\frac{2\varepsilon}{N}\right)^2+
        \left(s+\frac{\alpha}{{{N^{1/4}}}}\right)^2}
  = \sqrt{(r^2+s^2)+\frac{2}{{{N^{1/4}}}}\beta},
  \qquad &&\text{(\eqref{eq:length-1bis}, \eqref{eq:fp-proj})},\\
 \end{alignat*}
 where
\[
    \beta := s\alpha +     \left( \frac{\alpha^2}{2{{N^{1/4}}}}
 +\frac{2 r\varepsilon}{N^{3/4}}+ 
               \frac{ 2\varepsilon^2}{N^{7/4}}\right)
          \leq 
        s\alpha + 
       \frac{\alpha^2}{2}+2  r \varepsilon+
                2 \varepsilon^2  
           =  c_3\sqrt{r^2+s^2}
\]
 and $c_3$ is the constant as in \eqref{eq:c123}.
 Hence by the inequality $\sqrt{1+x}\leq 1+(x/2)$, 
 \begin{align*}
  |F_j(p)|&\leq \sqrt{r^2+s^2}\sqrt{1+\frac{2}{{{N^{1/4}}}}
               \frac{\beta}{r^2+s^2}}
          \leq 
             \sqrt{r^2+s^2}\left(1+\frac{1}{{{N^{1/4}}}}
            \frac{\beta}{r^2+s^2}\right)\\
          &\leq 
            \sqrt{r^2+s^2}\left(1+\frac{1}{{{N^{1/4}}}}
            \frac{c_3}{\sqrt{r^2+s^2}}\right)
          \leq \sqrt{r^2+s^2}+\frac{c_3}{{{N^{1/4}}}}
 \end{align*}
 holds,
 which is the conclusion.
\end{proof}
\begin{corollary}\label{cor:d3-case-1}
 Under the assumption of Lemma~\ref{lem:d3-case-1}, 
 we have
 \[
    |Y(p)|=|F_{2N}(p)|\leq \sqrt{r^2+s^2}+\varepsilon
    \qquad \text{for $p\in(\partial\DD\cap\varpi_j)$.}
 \]
\end{corollary}
\begin{proof}
\allowdisplaybreaks{
\begin{alignat*}{2}
   |F_{2N}(p)|&\leq 
   |F_j(p)|+|F_{2N}(p)-F_j(p)|
      \leq 
   |F_j(p)|+\frac{2\varepsilon}{N}
   \qquad &&\text{(by \eqref{eq:length-2} in Lemma~\ref{lem:length})}\\
   &\leq
   \sqrt{r^2+s^2}+\frac{c_3}{{{N^{1/4}}}}+\frac{2\varepsilon}{N}
  \qquad &&\text{(Lemma~\ref{lem:d3-case-1})}\\
   &= 
   \sqrt{r^2+s^2}+\frac{1}{{{N^{1/4}}}}
                    \left(c_3+\frac{2\varepsilon}{N^{3/4}}\right)\\
   &\leq 
   \sqrt{r^2+s^2}+\frac{1}{{{N^{1/4}}}}\left(c_3+2\varepsilon\right)
   \leq 
   \sqrt{r^2+s^2}+\varepsilon
 \qquad &&\text{(by \eqref{eq:N4})}.\qquad \qed
\end{alignat*}}%
\renewcommand{\qedsymbol}{\relax}
\end{proof}
\subsubsection*{Case 2: The case that $p\in\varpi_j$ 
       and  $|F_{j-1}(\zeta_j)|\leq 1/\sqrt{N}$}
\begin{lemma}\label{lem:d3-case-2}
 When $|F_{j-1}(\zeta_j)|\leq 1/\sqrt{N}$ and $p\in(\partial\DD\cap\varpi_j)$,
 \[
     |F_{2N}(p)|\leq \sqrt{r^2+s^2}+\varepsilon
 \]
 holds.
\end{lemma}
\begin{proof}
 Since $\bar p\in\partial\varpi_j$, 
\allowdisplaybreaks{
\begin{alignat*}{2}
 &|F_{2N}(p)|\leq |F_{2N}(p)-F_{2N}(\bar p)|+|F_{2N}(\bar p)|
            \leq \left(s+\frac{c_2}{{{N^{1/4}}}}\right)+|F_{2N}(\bar p)|
            && 
            \text{(by \eqref{eq:length-3})}\\
            &\leq \left(s+\frac{c_2}{{{N^{1/4}}}}\right)+
             |F_{2N}(\bar p)-F_{j}(\bar p)|+|F_{j}(\bar p)|\\
            &\leq 
             \left(s+\frac{c_2}{{{N^{1/4}}}}\right)+
              \frac{2\varepsilon}{N}
             +|F_{j}(\bar p)|
           \qquad&&\text{(by \eqref{eq:length-2})}\\
            &\leq 
                s + \frac{1}{{{N^{1/4}}}}\left(c_2+\frac{2\varepsilon}{N^{3/4}}\right)
             +|F_{j}(\bar p)-F_{j-1}(\bar p)|+|F_{j-1}(\bar p)|\\
            &\leq 
                s + \frac{1}{{{N^{1/4}}}}\left(c_2+\frac{2\varepsilon}{N^{3/4}}\right)
             +\frac{\varepsilon}{N^2}+|F_{j-1}(\bar p)|
           \qquad&&\text{(Lemma~\ref{lem:fj-est})}\\
            &\leq 
                s + \frac{1}{{{N^{1/4}}}}  \left(c_2+\frac{2\varepsilon}{N^{3/4}} +\frac{\varepsilon}{N^{7/4}}\right)
                  |F_{j-1}(\bar p)-F_{j-1}(\zeta_j)|+|F_{j-1}(\zeta_j)|\\
            &\leq 
                s + \frac{1}{{{N^{1/4}}}}  \left(  c_2+\frac{2\varepsilon}{N^{3/4}} +\frac{\varepsilon}{N^{7/4}}
               +   \frac{6\mu}{N^{3/4}} \right)+|F_{j-1}(\zeta_j)|
           \qquad&&\text{(Lemma~\ref{lem:diam-omega})}\\
            &
            \leq 
                s +  \frac{1}{{{N^{1/4}}}}  \left(  c_2+\frac{2\varepsilon}{N^{3/4}} +\frac{\varepsilon}{N^{7/4}}
               +   \frac{6\mu}{N^{3/4}} 
                +\frac{1}{N^{1/4}} \right)\\
      &= s+\frac{1+c_2+6\mu+3\varepsilon}{{{N^{1/4}}}}
             \leq s+\varepsilon
             \leq \sqrt{r^2+s^2}+\varepsilon
           \qquad&&\text{(by \eqref{eq:N4})}.
\end{alignat*}}%
Thus we have the conclusion.
\end{proof}

\subsubsection*{%
  Case 3:
  $p\in\overline{\D}_1\setminus(\varpi_1\cup\dots\cup\varpi_{2N})$
}\label{subsub:case-3}
The remaining case is that
 \begin{equation}\label{eq:p-remain}
    p\in \partial\DD \cap
         \bigl(
	  \overline{\D}_1\setminus(\varpi_1\cup\dots\cup\varpi_{2N})
	  \bigr).
 \end{equation}
\begin{lemma}\label{lem:d3-case-3}
 If $p$ satisfies \eqref{eq:p-remain},
 then
 \[
    |F_{2N}(p)|\leq \sqrt{r^2+s^2}+\varepsilon
 \]
 holds.
\end{lemma}  
\begin{proof}
 By the assumption \ref{key:ass:3}, $|X(p)|=|F_0(p)|\leq r$ holds.
 Then by Lemma~\ref{lem:d1} and \eqref{eq:N1},
 we have
 \[
    |F_{2N}(p)|\leq |F_0(p)|+|F_{2N}(p)-F_{0}(p)|\leq 
                    r + \frac{2\varepsilon}{N}
            \leq r+\varepsilon
           \leq \sqrt{r^2+s^2}+\varepsilon.\qquad \qed
 \]
\renewcommand{\qedsymbol}{\relax}
\end{proof}

Summing up, Corollary~\ref{cor:d3-case-1} and
Lemmas~\ref{lem:d3-case-2} and \ref{lem:d3-case-3}
implies \ref{key:concl:d3} of the Proposition~\ref{prop:main}.

\appendix
\section{Labyrinth}\label{app:labyrinth}
For the sake of completeness,  we recall Nadirashvili's labyrinth (for further details we refer the reader to
\cite{Nadi} or \cite{CR}).

For each number $j=0,1,2,\dots,2N^{2}$, we set
\begin{equation}\label{eq:rj}
    r_j := 1-\frac{j}{N^{3}}
    \qquad
    \left(
      r_0=1, r_1=1-\frac{1}{N^{3}},\dots,r_{2N^{2}}=1-\frac{2}{N}
    \right),
\end{equation}
and take a sequence of domains
\begin{equation}\label{eq:dj}
    \D_{r_j}=\{z\in\C\,;\,|z|<r_j\}
    \qquad (j=0,\dots,2N^{2}).
\end{equation}
Since $\{r_j\}$ is decreasing in $j$, it holds that
\[
   \D_1=\D_{r_0}\supset \D_{r_1}\supset \dots\supset
   \D_{r_{2N^{2}}}
   =\D_{1-\frac{2}{N}}.
\]
We denote the boundaries of $\D_{r_j}$ by 
\begin{equation}\label{eq:boundary-rj}
   S_{r_j}=\partial \D_{r_j}
               =\{z\in\C\,;\,|z|=r_j\}.
\end{equation}
We set 
\begin{align}\label{eq:all-annulus}
 \A &:= \overline{\D}_1\setminus
 \D_{r_{2N^{2}}}=\overline{\D}_1\setminus \D_{1-\frac{2}{N}}
\end{align}
and 
\begin{align*}
  A &:= \bigcup_{j=0}^{N^{2}-1} 
         \bigl(\D_{r_{2j}}\setminus \D_{r_{2j+1}}\bigr)
     =(\D_{r_0}\setminus \D_{r_1})\cup 
      (\D_{r_2}\setminus \D_{r_3})\cup \dots \cup
      (\D_{r_{2N^2-2}}\setminus \D_{r_{2N^{2}-1}}),
\nonumber\\
  \widetilde A &:= \bigcup_{j=0}^{N^{2}-1} 
         \bigl(\D_{r_{2j+1}}\setminus \D_{r_{2j+2}}\bigr)
     =(\D_{r_1}\setminus \D_{r_2})\cup 
      (\D_{r_3}\setminus \D_{r_4})\cup \dots \cup
      (\D_{r_{2N^2-1}}\setminus \D_{r_{2N^{2}}}).
\nonumber
\end{align*}
\begin{figure}
\begin{center}
\begin{tabular}{c@{\hspace{1cm}}c}
\includegraphics[width=0.4\textwidth]{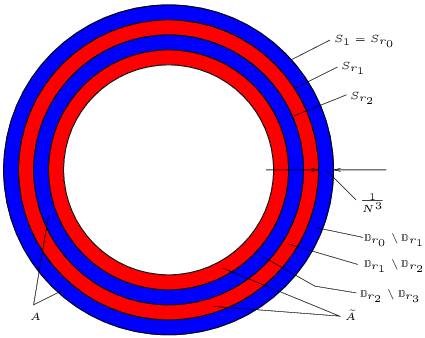}&
\includegraphics[width=0.4\textwidth]{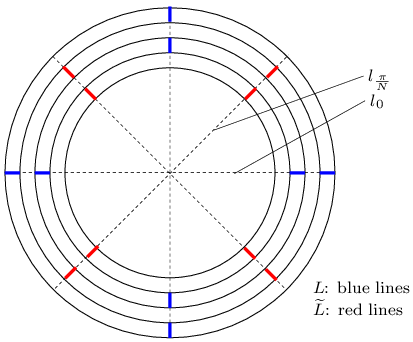}\\
\includegraphics[width=0.4\textwidth]{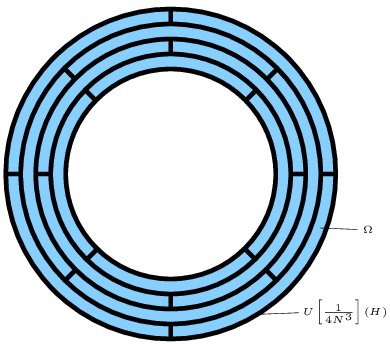}&
\includegraphics[width=0.4\textwidth]{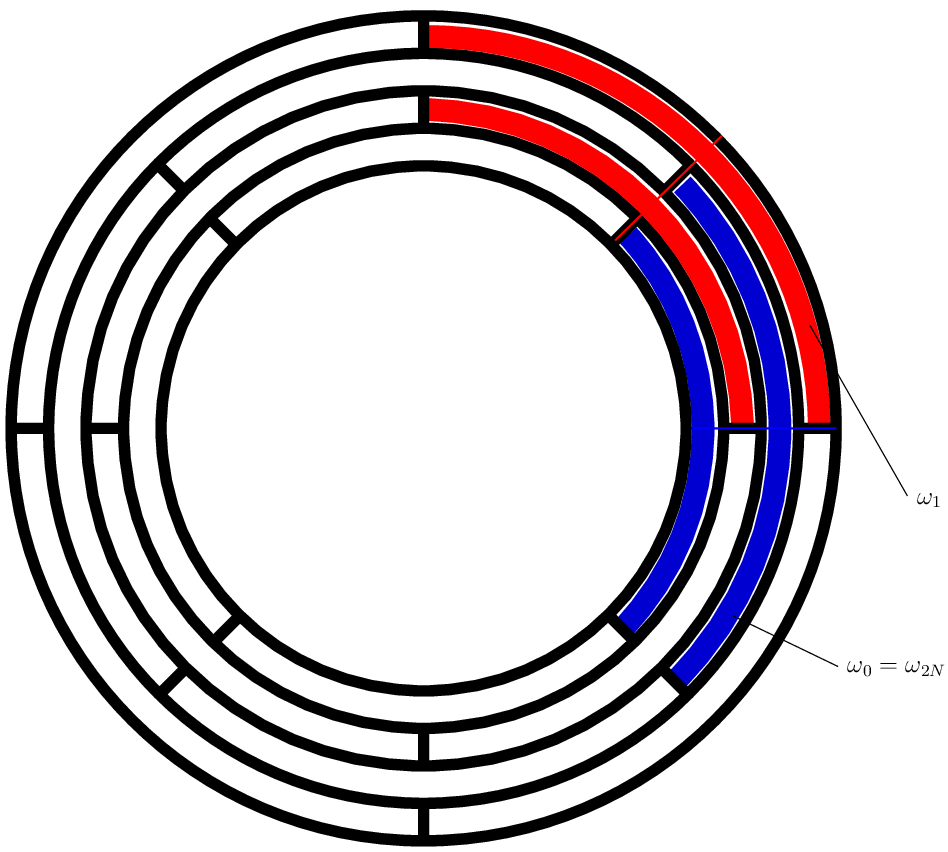}\\
\end{tabular}

\end{center}
\caption{The labyrinth}\label{fig:labyrinth}
\end{figure}
\begin{figure}
\begin{center}
 \includegraphics[width=0.4\textwidth]{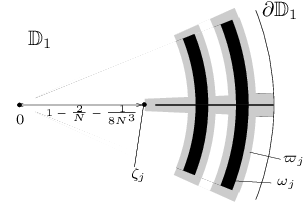}
\end{center}
\caption{The base point $\zeta_j$}
\label{fig:basepoint}
\end{figure}
Next, let 
\begin{equation}\label{eq:argument}
 L := \left(\bigcup _{j=0}^{N-1}l_{\frac{2j\pi}{N}}\right)\cap A,\qquad
 \widetilde L := 
     \left(\bigcup
      _{j=0}^{N-1}l_{\frac{(2j+1)\pi}{N}}\right)\cap\widetilde A,
\end{equation}
where $l_{t}:=\{re^{\imag t}\,;\,r\geq 0\}$,
and set 
\begin{equation}\label{eq:H}
 H:= L \cup \widetilde L \cup S\qquad
 \left( S = \bigcup_{j=0}^{2N^{2}} \partial \D_{r_j}
    = \bigcup_{j=0}^{2N^{2}} S_{r_j}
 \right).
\end{equation}
We define
\begin{equation}\label{eq:Omega}
 \Omega = \A \setminus U\left[\frac{1}{4N^{3}}\right](H),
\end{equation}
where $U[\varepsilon](B)$  denotes the $\varepsilon$-neighborhood of the subset
$B\subset \C$ (in the Euclidean distance).
Note that each connected component of $\Omega$ has the width 
$1/(2N^{3})$.

For each number $j=1,\dots,2N$, we set 
\begin{equation}\label{eq:omega}
\begin{aligned}
 \omega_j&:= 
   \bigl(l_{\frac{j\pi}{N}}\cap \A\bigr)
    \cup 
   \bigl(\text{the connected components of $\Omega$
          intersecting with $l_{\frac{j\pi}{N}}$}\bigr)\\
 \varpi_j&:=
   U\left[\frac{1}{8N^{3}}\right](\omega_j)
     =\left(\text{the $\frac{1}{8N^{3}}$-neighborhood of $\omega_j$}\right).
\end{aligned}
\end{equation}
Finally we denote by $\zeta_j$ the ``base point'' of $\varpi_j$:
\begin{equation}\label{eq:zeta}
   \zeta_j :=
    \left(1-\frac{2}{N}-\frac{1}{8N^{3}}\right)e^{\imag\pi j/N}
    \in\partial\varpi_j
    \qquad (j=1,\dots,2N)
\end{equation}
(Figure~\ref{fig:basepoint}).

By definition, we have
\begin{lemma}\label{lem:omega}
 \begin{enumerate}
  \item\label{item:omega:1}  
       For each $j=1,\dots,2N$,
       both $\omega_j$ and $\overline{\D}_1\setminus\varpi_j$
       are disjoint compact subsets of $\C$ such that
       $\C\setminus\bigl(\omega_j\cup(\overline{\D}_1\setminus\varpi_j)
       \bigr)$
       is connected.
  \item\label{item:omega:2}
       It holds that
       \[
	   \overline\D_1\setminus\D_{1-\frac{2}{N}-\frac{1}{8N^{3}}}
           \supset \varpi_{1}\cup\dots\cup\varpi_{2N}.
       \]
 \end{enumerate}
\end{lemma}
\begin{lemma}\label{lem:path}
 Let $j\in\{1,\dots,2N\}$.
 Then for each $p\in\overline{\D}_1\setminus\varpi_j$, there exists a path $\gamma$
 in $\overline{\D}_1\setminus\varpi_j$ joining $0$ and $p$
 whose length {\rm(}with respect to the Euclidean 
 metric of $\C${\rm)} is not greater than 
 $1+\pi/N$.
\end{lemma}
\begin{proof}
 By a rotation and a reflection on $\C=\R^2$,
 we assume $j=2N$ and 
 $p=re^{\imag\theta}$ ($0\leq r\leq 1$, $0\leq \theta\leq \pi$)
 without loss of generality.

 If $\pi/N<\theta<\pi$, 
 the line segment $\gamma$ joining $0$ and $p$ does not intersect 
 with $\varpi_{2N}$.
 Then $\gamma$ is the desired path.

 Otherwise, both the line segment $\gamma_1$ 
 joining $0$ and $p_0:=re^{\imag\pi/N}$ and 
 the circular arc $\gamma_2$ joining $p_0$ and $p$ centered at $0$
 do not intersect with $\varpi_{2N}$.
 Then the path $\gamma:=\gamma_1\cup \gamma_2$ is the desired one.
\end{proof}
\begin{lemma}\label{lem:omega-path}
 Let $j\in \{1,\dots,2N\}$.
 Then for each $p\in \overline{\varpi}_j$, there exists a path $\gamma$ in
 $\overline{\varpi}_j$
 joining the base point $\zeta_j$ and $p$
 whose length {\rm(}with respect to the Euclidean 
 metric of $\C${\rm)} is not greater than 
 $6/N$.
\end{lemma}
\begin{proof}
 We write $p=re^{\imag\theta}\in \overline{\varpi}_j$, where 
 \[
     1-\frac{2}{N}-\frac{1}{8N^{3}}\leq r\leq 1,\qquad
     \frac{\pi (j-1)}{N}\leq \theta\leq \frac{\pi (j+1)}{N}.
 \]
 Then the line segment $\gamma_1$
 joining $\zeta_j$ and $p_1:=re^{\imag\pi j/N}$ 
 lies in $\overline{\varpi}_j$, and its Euclidean length
 does not exceed  $\frac{2}{N}+\frac{1}{8N^{3}}$.
 On the other hand, the length of the circular arc $\gamma_2$
 centered at the origin
 joining  $p_1$ and $p$ does not exceed $\pi/N$.
 Then the path $\gamma=\gamma_1\cup\gamma_2$ joins $\zeta_j$
 and $p$ in $\overline{\varpi}_j$, whose length
 does not exceed 
 \[
    \frac{2}{N}+\frac{1}{8N^{3}}+\frac{\pi}{N}
    =\frac{1}{N}
         \left(2+\frac{1}{8N^{2}}+\pi\right)
    \leq \frac{1}{N}
         \left(2+\frac{1}{8}+\pi\right)\leq \frac{6}{N}.
 \]
 Hence we have the conclusion.
\end{proof}
\begin{lemma}\label{lem:length-estimate}
 Assume $N\geq 4$, and let 
 $\Omega\subset\D_1$ be the set as in  \eqref{eq:Omega}.
 Note that
 \[
    \Omega \subset \omega_1\cup\dots\cup \omega_{2N}.
 \]
 Consider a Riemannian metric $ds^2=\lambda^2\,|dz|^2$ on 
 $\overline{\D}_1$
 such that
 \[
    \begin{cases}
       \lambda \geq 1 \qquad & (\text{on $\overline{\D}_1$})\\
       \lambda \geq N^{3} \qquad & (\text{on $\Omega$}).
    \end{cases}
 \]
 Then for an arbitrary path $\sigma$ in $\overline{\D}_1$
 joining $0$ and $\partial\D_1$, it holds that
 $\int_{\sigma} ds \geq N$.
\end{lemma}
\begin{proof}  For $j=0,\ldots,N^2-1$, let $\gamma_j$ be a subarc 
 of $\sigma$ joining $\partial\D_{r_{2j}}$ and $\partial\D_{r_{2j+2}}$  contained in 
 $\overline{\D}_{r_{2j}}\setminus \D_{r_{2j+2}}$. 
 It suffices to prove that $\Length_{ds^2}(\gamma_j) \geq \frac{1}{N}$. In this case, 
 since the path $\sigma$ contains at least $N^2$ such paths,
 we have
$$ \Length_{ds^2}(\sigma) =\int_{\sigma}ds \geq N^2\cdot \frac{1}{N}=N, $$
In order to prove  that $\Length_{ds^2}(\gamma_j) \geq \frac{1}{N}$, we distinguish 
two cases. First we assume that $\Length_{\C}(\gamma_j) \geq \frac{1}{N}$. In this case
by the assumption $\lambda \geq 1$ we have
$$\Length_{ds^2}(\gamma_j) =\int_{\gamma_j}ds = \int_{\gamma_j}\lambda(z)\,|dz| \geq\int_{\gamma_j}|dz| \geq \frac{1}{N}$$ 
On the contrary, if $\Length_{\C}(\gamma_j) < \frac{1}{N}$ it is not difficult to see that $\gamma_j$ must be 
contained in a wedge of $\overline{\D}_1$ of angle bounded by $\frac{\pi}{N}-\frac{2}{N^2}$. Taking into account
the shape of the labyrinth, this implies that $\gamma_j$ crosses a connected component of $\Omega$
 transversely,
and therefore the Euclidean length of $\gamma_j\cap\Omega$ is greater than
 $1/(2N^{3})$.
 Hence by the assumption,
\[
   \Length_{ds^2}(\gamma_j)= \int_{\gamma_j}\,ds 
    \geq \int_{\gamma_j\cap\Omega}\,ds
    =\int_{\gamma_j\cap\Omega}\lambda\,|dz|
    \geq N^{3}\cdot\frac{1}{2N^{3}}=\frac{1}{2}> \frac{1}{N}.\qed
\]
\renewcommand{\qed}{\relax}
\end{proof}

\end{document}